\documentclass[10pt,leqno,a4paper]{article}
\usepackage{amsmath,amsfonts,amscd,amsthm}
\usepackage{graphicx}
\usepackage[dvipsnames]{xcolor}
\usepackage{cite}
\usepackage{tikz}
\usepackage{ulem}
\usepackage{hyperref}

\newcommand{\R}{\mathbb{R}}

\newcommand{\pt}{\partial_t}
\newcommand{\po}{\partial\Omega}

\newcommand{\bw}{{\bf w}}

\newcommand{\bv}{{\bf v}}

\newcommand{\bH}{{\bf H}}
\newcommand{\bL}{{\bf L}}

\newcommand{\ve}{\varepsilon}

\newcommand{\g}{{\bf g}}
\newcommand{\f}{{\bf f}}

\newcommand{\bu}{{\bf u}}

\newcommand{\bue}{{\bu}_\ve}
\newcommand{\buet}{{\bu}_{\ve t}}
\newcommand{\bues}{{\bu}_{\ve s}}
\newcommand{\buett}{{\bu}_{\ve tt}}

\newcommand{\hbue}{\hat{\bu}_\ve}
\newcommand{\hbuet}{\hat{\bu}_{\ve t}}

\newcommand{\pe}{{p}_\ve}

\newcommand{\tb}{\tilde{b}}
\newcommand{\hf}{\hat{\bf f}}

\newcommand{\bueh}{{\bu}_{\ve h}}
\newcommand{\peh}{{p}_{\ve h}}
\newcommand{\bueht}{{\bu}_{\ve{h t}}}
\newcommand{\buehtt}{{\bu}_{\ve{h tt}}}
\newcommand{\buehss}{{\bu}_{\ve{h ss}}}
\newcommand{\buehs}{{\bu}_{\ve{h s}}}
\newcommand{\hbueh}{\hat{\bu}_{\ve h}}
\newcommand{\hbueht}{\hat{\bu}_{\ve ht}}

\newcommand{\btveh}{\tilde{\bv}_{\ve h}}
\newcommand{\btveht}{\tilde{\bv}_{\ve h t}}

\newcommand{\bU}{{\bf U}}

\newcommand{\bUe}{{\bf U}_{\ve}}

\newcommand{\Pe}{P_{\ve}}

\newcommand{\bphi}{\mbox{\boldmath $\phi$}}
\newcommand{\bpsi}{\mbox{\boldmath $\psi$}}
\newcommand{\bxi}{\mbox{\boldmath $\xi$}}

\newcommand{\bta}{\mbox{\boldmath $\eta$}}

\newcommand{\bzeta}{\mbox{\boldmath $\zeta$}}
\newcommand{\btheta}{\mbox{\boldmath $\theta$}}

\newcommand{\eve}{{\bf e}_\ve}
\newcommand{\evet}{{\bf e}_{\ve t}}

\DeclareSymbolFont{myletters}{OML}{ztmcm}{m}{it}
\DeclareMathSymbol{\uplambda}{\mathord}{myletters}{"15}

\newcommand{\Seh}{{S}^\ve_h}
\DeclareMathOperator*{\esssup}{ess\,sup}

\begin{document}

\newcommand{\se}{\setcounter{equation}{0}}
\def\theequation{\thesection.\arabic{equation}}

\newtheorem{theorem}{Theorem}[section]
\newtheorem{cdf}{Corollary}[section]
\newtheorem{lemma}{Lemma}[section]
\newtheorem{remark}{Remark}[section]
\newtheorem{example}{Example}[section]
\newtheorem{proposition}{Proposition}[section]

\def\cydot{\leavevmode\raise.4ex\hbox{.}}

\title
{ \large\bf Optimal error estimates of the penalty finite element method for the unsteady Navier-Stokes equations with nonsmooth initial data}

\author{Bikram Bir\thanks{Department of Mathematics, Indian Institute of Technology Bombay, Powai, Mumbai-400076, India. Email: bikram@math.iitb.ac.in}, ~~
Deepjyoti Goswami\thanks{Department of Mathematical Sciences, Tezpur University, Tezpur, Sonitpur, Assam-784028, India. Email: deepjyoti@tezu.ernet.in}, ~~and
Amiya K. Pani\thanks{Department of Mathematics, BITS-Pilani, KK Birla Goa Campus, NH 17 B, Zuarinagar Goa-403726, India. Email: amiyap@goa.bits-pilani.ac.in,  akp@math.iitb.ac.in}}

\date{}
\maketitle

\begin{abstract}
In this paper, both semidiscrete and fully discrete finite element methods are analyzed for the penalized two-dimensional unsteady Navier-Stokes equations with nonsmooth initial data. First order backward Euler method is applied for the time discretization, whereas conforming finite element method is used for the spatial discretization. 
Optimal $L^2$ error estimates for the semidiscrete as well as the fully discrete approximations of the velocity and of the pressure are derived for realistically assumed conditions on the data. 
The main ingredient in the proof is the appropriate exploitation of the inverse of the penalized Stokes operator, negative norm estimates and time weighted estimates.
Two numerical examples one in 2D and one in 3D are presented whose results are conforming our theoretical findings. Finally, computational experiments on benchmark problem: one on lid driven cavity problem and other on flow around a cylinder with low viscosity are discussed.
\end{abstract}

\vspace{1em} 


\section{Introduction}
Let $\Omega$ be a bounded convex polygonal domain in $\mathbb{R}^2$ with boundary $\partial \Omega$. Now
consider the following  incompressible Navier-Stokes system in a space-time domain $\Omega\times (0,\infty)$
\begin{eqnarray}\label{om04}%------NSE
~~\frac {\partial \bu}{\partial t}+\bu\cdot\nabla\bu-\nu\Delta\bu +\nabla p=\f \quad\mbox {in}~ \Omega,~t>0,
\end{eqnarray}
with incompressibility condition
\begin{eqnarray}\label{ic04}
 \nabla \cdot \bu=0 \quad\mbox {on}~\Omega,~t>0,
\end{eqnarray}
and initial with boundary conditions
\begin{eqnarray}\label{ibc04}
 \bu(x,0)=\bu_0 \quad\mbox {in}~\Omega,\quad\bu=0 \quad\mbox {on}~\partial\Omega,~t\ge 0.
\end{eqnarray}
Here, $\bu$ denotes the velocity vector, $p$ represents the pressure of the fluid and $\nu>0$ is the kinematic coefficient of viscosity. Further, the forcing term $\f$ and the initial velocity $\bu_0$ are given functions in their respective domains of definition. 

The time-dependent Navier-Stokes equations (NSEs) for the incompressible flow has always been a major challenge in computational PDEs. The main difficulty in computation is that, 
at each time step, the velocity $\bu$ and the pressure $p$ are coupled together by the incompressibility condition, $div~ u = 0$. A common way to tackle this difficulty is to address the imposition of the incompressibility condition in an appropriate way so as to obtain a psedo-compressible system. In this regard, methods, which come to our mind, are the penalty method, the artificial compressibility method, the pressure stabilized method, the pressure correction method, the projection method, etc. (see, \cite{A20, AL14, AS15, BGT88, DTW13, FGN19, H05, HL10, H14, HR19, HHF13, HLY06, LA12, LA13, LL10, QZMX17, S92, S92_1, S95,T68, T84, ZKO16} and references, therein).

In the present paper, a completely discrete penalty finite element method to the NSEs is applied  to circumvent this difficulty. 
It is the simplest and effective finite element implementation to handle the incompressibility. This method is often employed  in order to decouple the pressure equation from the system of nonlinear algebraic equations in velocity which is obtained from finite element (or finite difference) discretizations of the Navier-Stokes equations at each time level.
The basic idea of the proposed method is to add a pressure term to the continuity equation.
The resulting penalized system in our case is to approximate the solution $(\bu,p)$ of (\ref{om04})-(\ref{ibc04}) by $(\bue,\pe)$ satisfying
\begin{equation}\left\{\begin{array}{l}\label{pom}
 \frac{\partial\bue}{\partial t}+\bue\cdot\nabla\bue-\nu\Delta\bue+\frac{1}{2}(\nabla\cdot\bue)\bue +\nabla\pe=\f(x,t) \quad\mbox {in}~ \Omega,~t>0,  \\
 \nu\nabla\cdot\bue+\ve\pe=0 \quad\mbox {on}~\Omega,~t>0,\quad~~
 \bue|_{t=0}=\bu_{\ve 0} \quad\mbox{in}~\Omega, \quad~~ \bue=0 \quad\mbox {on}~\partial\Omega,~t\ge 0.
\end{array}\right.
\end{equation}
Here $\ve>0$ is the penalty parameter. 
One advantage of this approach is that we can eliminate pressure $\pe$ from (\ref{pom}) to obtain an equation $\bue$ only.
We note here that it is standard to add the term $\frac{1}{2}(\nabla\cdot\bue) \bue$ to the non-linear term, which was introduced by Temam \cite{T68}, that ensures the dissipativity of the system (\ref{pom}).

The penalty method first appeared in the work of Courant \cite{C43} for a membrane problem, in the framework of the calculus of variation. Since then, it has been considerably developed over the time. Besides its applications to constrained variational problems and variational inequalities, the penalty method was proved to be quite successful
for numerical computations in continuum fluid and solid mechanics. Its application to NSEs has been initiated in the late 1960's by Temam \cite{T68}.
Subsequently, many publications have been devoted to the study of the penalty method for the steady Stokes and NSEs, as well as for the unsteady NSEs, in continuous, semidiscrete and fully discrete cases, see, for example,  \cite{BGT88, GR80, H05, HL10, S95} and references, therein. Recent results on the penalty method can be categorised as follows: two-grid penalty method \cite{AS15, HHF13, AL14}, iterative penalty method \cite{DTW13, H14}, methods based on different boundary conditions such as nonlinear slip boundary conditions \cite{AL14, DTW13, LA12}, friction boundary conditions \cite{LA13, QZMX17}, slip boundary conditions \cite{ZKO16}, etc.  
Moreover, penalty method is used, very recently, for the stochastic $2$-D incompressible NSEs \cite{HR19}, and for the incompressible NSEs with variable density \cite{A20}.

This article deals with optimal $L^2$-error estimates for the velocity and for the pressure term in the semidiscrete as well as in the fully discrete case when penalty finite element method is applied to NSE (\ref{om04})-(\ref{ibc04}) with nonsmooth initial data. 
Earlier, Shen \cite{S95} has applied the backward Euler method to discretize only in time and has derived optimal error estimates with respect to the penalty parameter and also with respect to the discretizing parameter in time.
Later on, He \cite{H05} has extended it to space discretization  and the following error estimate has been established for the conforming fully discrete finite element method for all $t_n\in [0,T],  T>0,$
\vspace*{-2mm}
\begin{align}\label{shen}
\tau(t_n)\|\bu(t_n)-\bueh^n\|_{H^1} + \Big( k\sum_{m=0}^n \tau^2(t_m)\|p(t_m)-p_{\ve h}^n\|_{L^2}^2	\Big)^{\frac{1}{2}} \le C(\ve + h+ k),
\end{align}
where $(\bu(t_n),p(t_n))$ and $(\bueh^n,\peh^n)$ are the solutions of the NSEs and its fully discrete penalized system, respectively. Here, $C$ is the positive constant, $h$ is the mesh size, $k$ is the time step, $t_n=nk, 0\le n \le N=T/k, \tau(t_n)=\min\{t_n,1\}.$  In both these paper, smallness condition on $k$ is assumed like $Ck<1$, where $C$ depends on $\nu^{-1}$, Sobolev constants, etc.
Subsequently, Lu and Lin \cite{LL10} have discussed optimal error estimates for the semidiscrete problem under smooth initial data.
But to the best of our knowledge, optimal $L^2$  error estimates for the velocity and the pressure in the semidiscrete and fully discrete cases of the penalized unsteady NSEs with nonsmooth initial data have not been obtained in the literature although these results have been established numerically on various occasions.
The purpose of this paper is to fill this gap. We extend the work of Shen \cite{S95} and He \cite{H05}, and obtain the optimal error estimates for the velocity and the pressure in $L^2$-norm.

Our analysis is based under realistically assumed conditions on the initial data $\bu_0$ in $\bH_0^1$. We take into account the lack of regularity endured by the solutions of the Navier-Stokes system at the initial time $t=0$. Assuming otherwise requires the data to satisfy some non-local compatibility conditions, which are not natural and difficult to verify in practice, see \cite{HR82}.
We note that in \cite{S95}, the penalized error estimates have been obtained for $\bu_0$ in $\bH_0^1$, but the time discrete error estimates have been obtained under additional assumption of $\bu_0$ in $\bH^2\cap\bH_0^1$. 
Also in \cite{H05}, the results have been obtained for the smooth initial data ($\bu_0$ in $\bH^2\cap\bH_0^1$). We take a more realistic approach and consider nonsmooth initial data, that is, $\bu_0$ in $\bH_0^1$, which poses more serious difficulties, mainly in both semidiscrete and fully discrete analyses. 

\noindent
The main results of this article consist of the following:
\begin{itemize}  
\item uniform in time bounds for the semidiscrete as well as fully discrete solution are established with no additional smallness assumptions on time step $k$.
\item optimal error estimates for the finite element approximation of the penalized velocity and pressure  are derived,
\begin{align*}
 \|\bue(t)-\bueh(t)\| + h \big(\|\nabla(\bue(t)-\bueh(t))\| + \|\pe(t)-\peh(t)\|\big) \leq C  h^{m+1}t^{-\frac{m}{2}},
% \|\bue(t)-\bueh(t)\| \le C  h^{m+1}t^{-\frac{m}{2}}, \\
% \|\nabla(\bue(t)-\bueh(t))\| \le C  h^{m}t^{-\frac{m}{2}},\\
% \|\pe(t)-\peh(t)\| \leq C  h^{m}t^{-\frac{m}{2}} .
\end{align*}
where $(\bue,\pe)$ and $(\bueh,\peh)$ are the solution of (\ref{pom}) and it's semidiscrete system, respectively. 
\item  optimal error bounds for the velocity and the pressure terms, for the fully discrete penalized system, are derived which are of the form:
\begin{align*}
 \|\bu(t_n)-\bUe^n\| & \le C\Big((\ve+k)t^{-\frac{1}{2}} + h^{m+1}t^{-\frac{m}{2}}\Big), \\
 \|\nabla(\bu(t_n)-\bUe^n)\|+\|p(t_n)-P_{\ve}^n\| & \le C\Big((\ve+k)t^{-1} + h^{m}t^{-\frac{m}{2}}\Big),
% \|p(t_n)-P_{\ve}^n\| \leq C\Big((\ve+k)t^{-1} + h^{m}t^{-\frac{m}{2}}\Big),
\end{align*}
where $(\bu,p)$ and $(\bUe^n,\Pe^n)$ are the solution of (\ref{om04})-(\ref{ibc04}) and the fully discrete system of (\ref{pom}), respectively.
\item since constants involved in the error estimates of both semidiscrete and fully discrete schemes depend exponentially on time, using uniqueness assumption, the error estimates are shown to be valid uniformly in time.
\item a couple of numerical examples are discussed to verify the theoretical findings and another couple of examples on benchmark problems with small viscosity are presented.
\end{itemize}

The remaining part of this paper is arranged as follows: Notations, assumptions and a couple of standard results are stated in the first part of the Section 2, whereas in the second part, we briefly look at the penalty method.
The semidiscrete error analysis is carried out in Section 3 and in Section 4, backward Euler method is applied to the penalized system. Finally, in Section 5, some numerical examples are given which validate our theoretical findings.

\section{Preliminaries}
\se

For our subsequent use, we denote by bold face letters, the $\R^2$-valued function space such as $\bH_0^1 = [H_0^1(\Omega)]^2$, $ \bL^2 = [L^2(\Omega)]^2$ and $\bH^m=[H^m(\Omega)]^2$.
We denote by $\|\cdot\|_m$ the usual norm of the Sobolev space $\bH^m$, and $(\cdot,\cdot)$ and $\|\cdot\|$ represent the inner product and norm on $L^2$ or $\bL^2$, respectively. 
The space $\bH^1_0$ is equipped with the norm
\[ \|\nabla\bv\|= \big(\sum_{i,j=1}^{2}(\partial_j v_i, \partial_j
 v_i)\big)^{\frac{1}{2}}=\big(\sum_{i=1}^{2}(\nabla v_i, \nabla v_i)\big)^{\frac{1}{2}}. \]
Let $H^m/\mathbb{R}$ be the quotient space of equivalent classes of functions in $\bH^m$ differ by constant with norm $\|\phi\|_{H^m/\mathbb{R}}=\inf_{c\in \mathbb{R}}\|\phi+c\|_m.$
For $m=0$, it is denoted by $L^2/\mathbb{R}$. 
For any Banach space $X$, let $L^p(0, T; X)$ denote the space of measurable $X$-valued functions $\bphi$ on  $ (0,T) $ such that
\[ \int_0^T \|\bphi (t)\|^p_X~dt <\infty, ~~\mbox {if}~1 \le p < \infty, \quad\mbox{and}\quad \esssup_{0<t<T} \|\bphi (t)\|_X <\infty, ~~\mbox {if}~p=\infty. \]
The dual space of $H^m(\Omega)$, denoted by $H^{-m}(\Omega)$, is defined as the completion of $C^\infty(\bar{\Omega})$ with respect to the norm 
$$\|\phi\|_{-m}:=\sup \bigg\{\frac{(\phi,\psi)}{\|\psi\|_m}: \psi \in H^m(\Omega), \|\psi\|_m\neq 0 \bigg\}.$$
Through out this paper, we make the following assumption:\\
(${\bf A1}$) For ${\bf g} \in \bH^{m-1}$ with $m\geq 1$, let the unique pair of solution $(\bv, q)\in\bH_0^1\times L^2/\R$ for the steady state Stokes problem
\begin{align*}
 & -\Delta\bv + \nabla q = {\bf g}, \quad \nabla \cdot\bv = 0,~~\mbox {in}~\Omega, \quad\bv|_{\partial\Omega}=0,
\end{align*}
satisfies $(\bv, q)\in\bH^{m+1}\times H^{m}/\R$ and the regularity result \cite{HR90}: $\|\bv\|_{m+1} + \|q\|_{H^m/\R} \le C\|{\bf g}\|_{m-1}.$\\
Assumption ($\textbf{A1}$) simply talk about the regularity of the boundary $\partial\Omega$, that is $\partial\Omega\in C^m$.
We note here that ({\bf A1}) implies
\begin{equation*}
\|\bv\|_{m+1}\le C\|\Delta\bv\|_{(m+1)/2},~~\forall\bv\in\bH_0^1\cap\bH^{m+1} \quad \mbox{and} \quad
\|\bv\|_{m-1}\le \lambda_1^{-\frac{1}{2}}\|\bv\|_m,~~\forall\bv\in\bH_0^1\cap\bH^m,
\end{equation*}
where $\lambda_1>0$ to be the least eigenvalue of the Stokes operator.

\noindent 
We present below a couple of lemmas for subsequent use.
\begin{lemma} [uniform Gronwall's Lemma \cite{BGT88}]\label{ugl}
Let $g,h,y$ be three locally integrable non-negative functions  on the time interval $[0,\infty)$.
Assume that $y$ is absolutely continuous and
 $$ \frac{dy}{dt} \le gy + h, \quad \forall t\ge 0,$$
and 
$$\int_t^{t+T} g(s)ds \le \alpha_1, \quad \int_t^{t+T} h(s)ds \le \alpha_2, \quad  \int_t^{t+T} y(s)ds \le \alpha_3, \quad \forall t\ge 0,$$
where $T, \alpha_1,\alpha_2, \alpha_3$ are positive constants. Then 
 $$ y(t+T)\le \Big(\frac{\alpha_3}{T}+\alpha_2\Big)\exp(\alpha_1),\quad \forall t\ge 0. $$
\end{lemma}
\begin{lemma}[discrete uniform Gronwall's Lemma \cite{S90}]\label{dugl}
Let $k$ and $\{y^i, g^i,h^i\}_{i\in\mathbb{N}}$ be non-negative numbers satisfying
\begin{equation*}%\label{dg1}
\frac{y^{i}-y^{i-1}}{k}  \le  g^{i-1} y^{i-1} + h^{i-1}, \quad \forall i\ge 1,
\end{equation*}
and there exist $a_1(r), a_2(r),a_3(r)$ depend on $t_r=rk$ such that
\begin{equation*}
k\sum_{i=i_0}^{i_0+r} g^i \le a_1(r),\quad
k\sum_{i=i_0}^{i_0+r} h^i \le a_2(r),\quad
k\sum_{i=i_0}^{i_0+r} y^i \le a_3(r), \quad \forall i_0\ge 1.
\end{equation*}
Then,
\begin{equation*}%\label{dg2}
y^n \le \Big\{\frac{a_3(r)}{t_r}+a_2(r)\Big\}\exp(a_1(r)), \quad \forall n\ge r+1.
\end{equation*}
\end{lemma}

We are now in a position to look at the variational formulation of the penalized system (\ref{pom}). The corresponding variational formulation of penalized NSEs is to find $(\bue(t),\pe(t)),~t>0$ in $\bH_0^1\times L^2$ satisfying
\begin{equation}\label{wfpom}
\left\{
\begin{array}{rl}
(\buet,\bphi)+\nu a(\bue,\bphi)+\tb(\bue,\bue,\bphi) -(\pe,\nabla\cdot\bphi) &=(\f,\bphi), \quad\forall~\bphi\in\bH_0^1, \\
 \nu(\nabla\cdot\bue,\chi)+\ve(\pe,\chi) &=0, \quad\forall~\chi\in L^2,
\end{array}\right.
\end{equation}
with $\bue(0)=\bu_{\ve 0}$. Here
$ a(\bv,\bw)= (\nabla\bv,\nabla\bw)$ and
\[ \tb(\bv,\bw,\bphi)= (\tilde{B}(\bv,\bw),\bphi),~~~ \mbox{where}~~ \tilde B(\bv,\bw) := (\bv\cdot\nabla)\bw+\frac{1}{2}(\nabla\cdot\bv)\bw. \]
We can easily check with the help of integration by parts that
$$ \tb(\bv,\bw,\bphi)=\frac{1}{2}\big\{b(\bv,\bw,\bphi)-b(\bv,\bphi,\bw)\big\}, \quad\forall ~\bv,\bw,\bphi\in\bH_0^1, $$ 
where $b(\bv,\bw,\bphi)= ((\bv\cdot\nabla)\bw,\bphi)$.
Hence, it follows that
\begin{equation}\label{tbp1}
\tb(\bv,\bw,\bw)=0,\quad \text{and}\quad \tb(\bv,\bw,\bphi)=-\tb(\bv,\bphi,\bw),\quad\forall~\bv,\bw,\bphi\in\bH_0^1.
\end{equation}
In order to omit the pressure term, we choose $\chi=\nabla\cdot\bphi$ in the second equation of (\ref{wfpom}). We then obtain 
\begin{equation}\label{wfpp}
(\buet,\bphi)+\nu a_{\ve}(\bue,\bphi)+\tb(\bue,\bue,\bphi)
=(\f,\bphi), \quad\forall~\bphi\in\bH_0^1,
\end{equation}
where $$  a_{\ve}(\bv,\bw)= a(\bv,\bw)+\frac{1}{\ve}(\nabla\cdot\bv,\nabla\cdot\bw), $$  with $\bue(0)=\bu_{\ve 0}.$ 
Setting $A_{\ve}\bv := -\Delta\bv-\frac{1}{\ve}\nabla (\nabla\cdot\bv)$, we rewrite the system in abstract form as
\begin{equation}\label{pomu}
 \buet +\nu A_{\ve}\bue+\tilde B(\bue,\bue)=\f.
\end{equation}

The operator $A_{\ve},$ which is associated with the penalty method, is a self-adjoint and  positive definite operator from $\bH^2\cap \bH_0^1$ onto $\bL^2.$
Similar to the Stokes operator, we can talk of various powers of $A_{\ve}$, namely, $A_{\ve}^r, ~r\in\mathbb{R}.$ For details, we refer to Temam \cite{BGT88} and Shen \cite{S95}.
It is observed in \cite{BGT88} that $\|A_{\ve}\bv\|$ is a norm on $\bH^2\cap\bH_0^1$ and is, in fact, equivalent to that of $\bH^2,$ i.e.,
\begin{equation}\label{aeeqv}%-----------------ae eqv norm
 \|A_{\ve}\bv\|\cong\|\bv\|_2,
\end{equation}
with constants depending on $\ve$. In \cite{BGT88}, one of the inequalities of (\ref{aeeqv}), that is, $\|\Delta\bv\|\le C_3\|A_{\ve}\bv\|$, for some positive constant $C_3$ with $\ve C_3\le 1$  is proved to be independent of $\ve$, which is crucial for our subsequent analysis.
We present below a Lemma, to support this. For a proof, we again refer to \cite[pp. 6]{BGT88} and \cite[pp. 388]{S95}.
\begin{lemma}\label{Aep}%------------------------------------------A_{\ve} prop.
 There exists a constant $c_0>0$ such that, for $\ve>0$ sufficiently small, the following 
estimates hold:
\begin{align*}
& \|\bv\|_{r} \le  c_0\|A^{\frac{r}{2}}_{\ve}\bv\|, \quad\forall~\bv\in \bH^{r}\cap \bH_0^1,~~1\le r\le m+1, \\
 & \|A_{\ve}^{-1}\bv\| \le  c_0\|\bv\|_{-2},\quad\forall~\bv\in\bH^{-2}.
\end{align*}
\end{lemma}
\noindent
Based on the second result, we have obtained another estimate, independent of $\ve,$ which we prove in the next Lemma.
\begin{lemma}\label{Aep1}%------------------------------------------A_{\ve} prop. 1
 There exists a constant $c_0>0$ such that, for $\ve>0$ sufficiently small, the following holds 
true:
$$ \|A_{\ve}^{-\frac{1}{2}}\bv\| \le c_0\|\bv\|_{-1}, \quad\forall~\bv\in\bH^{-1}. $$
\end{lemma}

\begin{proof}
With $\bw\in\bH_0^1,$ we use Lemma \ref{Aep} to find that
$$ (A_{\ve}^{-\frac{1}{2}}\bv, \bw)= (\bv, A_{\ve}^{-\frac{1}{2}}\bw) \le \|\bv\|_{-1}\| \nabla (A_{\ve}^{-\frac{1}{2}} 
\bw)\| \le c_0\|\bv\|_{-1}\|\bw\|, $$
and
$$ \|A_{\ve}^{-\frac{1}{2}}\bv\| = \sup_{0\neq \bw\in\bL^2} \frac{(A_{\ve}^{-\frac{1}{2}}\bv, 
\bw)}{\|\bw\|} \le c_0\|\bv\|_{-1}. $$
This completes the proof. Alternative way is to consider the following problem:
Let $\bw$ be a solution of
$$ A_{\ve}\bw=\bv,~~\bw|_{\po}=0. $$
Clearly $\|A_{\ve}\bw\|=\|\bv\|$ and
$$ \|A_{\ve}^{\frac{1}{2}}\bw\|^2=(A_{\ve}\bw,\bw)=(\bv,\bw)=(A_{\ve}^{-\frac{1}{2}}\bv,A_{\ve}^{\frac{1}{2}}
   \bw) \le \|A_{\ve}^{-\frac{1}{2}}\bv\|\|A_{\ve}^{\frac{1}{2}}\bw\| $$
and therefore
\begin{equation}\label{Aep101}
 \|A_{\ve}^{\frac{1}{2}}\bw\| \le \|A_{\ve}^{-\frac{1}{2}}\bv\|.
\end{equation}
Now using (\ref{Aep101}) and Lemma \ref{Aep}, we note that
\begin{align*}
 \|A_{\ve}^{-\frac{1}{2}}\bv\|^2 &=(A_{\ve}^{-1}\bv, \bv)=(A_{\ve}^{-1}\bv, A_{\ve}\bw) = 
(\bv,\bw) \le \|\bv\|_{-1}\|\bw\|_1 \\
& \le c_0\|\bv\|_{-1}\|\|A_{\ve}^{\frac{1}{2}}\bw\| \le c_0\|\bv\|_{-1}\| \|A_{\ve}^{-\frac{1}{2}}\bv\|.
\end{align*}
This completes the rest of the proof.
\end{proof}
\noindent  
We now consider the following assumptions on the given data for the penalized NSEs for our subsequent analysis.
\\
%%-------------------------------------------------------------------------A4
\noindent 
(${\bf A2}$).
 {\it The initial velocity $\bu_{\ve 0}$ and the external force $\f$ satisfy for positive
constant $M_0,$ and for $T$ with $0<T < \infty$} and for some integer $m\geq 1$ and $0\le l \le m$
$$ \bu_{\ve 0}\in\bH_0^1(\Omega),~ D^{l}_{t}\f \in L^{\infty}(0,T;{\bH}^{m-1})~~~\mbox{with}~~~
 \|A_{\ve}^{\frac{1}{2}}\bu_{\ve 0}\| \le M_0,~~{{\sup_{0<t<T} }}\Big\{\|D^{l}_{t}\f(t)\|_{m-1}\Big\} \le M_0. $$ 

\noindent 
Throughout, we shall use $C$ as a generic constant depending on the data: $ \Omega, \bu_{\ve 0}, \f, \nu, c_0$ and $T$, but not on mesh parameter $h$ and $k$.
Below, we take a quick glance at the {\it a priori} estimates of the penalized problem.
\begin{lemma}\label{pap}
 Assume $({\bf A1})$ and $({\bf A2})$ hold true and $0< \alpha <\nu \lambda_1/2c_0^2$. Then,
there exists constant $C>0,$ independent of $\ve,$ such that
\begin{eqnarray*}
 \|\bue(t)\|^2 + e^{-2\alpha t}\int_0^t e^{2\alpha s} \|A_{\ve}^{\frac{1}{2}}\bue(s)\|^2~ds &\le & C , \\
\tau^{m}(t)\|A_{\ve}^{\frac{m+1}{2}}\bue(t)\|^2 + e^{-2\alpha t}\int_0^t \sigma^{m}(s)\|A_{\ve}^{\frac{m+2}{2}} \bue(s)\|^2 ds &\le & C, \\
\tau^{m+1}(t)\|A_{\ve}^{\frac{m}{2}}\buet(t)\|^2+ e^{-2\alpha t}\int_0^t \sigma^{m+1}(s)\|A_{\ve}^{\frac{m+1}{2}} \bues(s)\|^2  &\le & C, 
\end{eqnarray*}
hold for $m\le 4$, where $\tau(t)=min\{t,1\}$ and $\sigma^{m}(t)=\tau^m(t) e^{2\alpha t}$.
\end{lemma}
\noindent 
The proof goes in a similar way as that of the proofs given in \cite[Lemma 2.1]{GD15}  and \cite[Proposition 3.2]{HR90}. We give a sketch in the appendix.\\
%
%%------------------------------------------PENALIZED PROBLEM: ERROR ESTIMATE
For the sake of completeness, we present below the optimal penalty error estimate.% 
\begin{theorem}\label{pperrest}
 Under the assumption of Lemma \ref{pap}, the following hold:
\begin{equation*} %\label{pperrest1}
 \sqrt{\tau(t)}\|(\bu-\bue)(t)\|+ \tau(t)\|\nabla(\bu-\bue)(t)\| + \Big(
 \int_0^t \tau^2(s)\|(p-\pe)(s)\|^2ds\Big)^{\frac{1}{2}} \le K(t) \ve,
\end{equation*}
and 
$$\tau(t) \|(p-\pe)(t)\| \leq K(t) \ve,$$
where $\tau(t)=\min\{t,1\}$ and $K(t)=Ce^{Ct}$, $C$ is the positive constant. 
Under the uniqueness condition:
\begin{equation}\label{unique}
\frac{2N}{\nu^2} \|\f\|_{L^\infty(0,\infty;\bH^{-1}(\Omega))}<1 ~\text{and}~ N= \sup_{\bu,\bv,\bw}\frac{b(\bu,\bv,\bw)}{\|\nabla\bu\|\|\nabla\bv\|\|\nabla\bw\|},
\end{equation}
the above estimates are uniform in time, that is, the constant $K(t)$ becomes $C$.
\end{theorem}
\begin{proof}
The proof of the first estimate is available in \cite[Theorem 4.1]{S95}. Using this, we can easily prove the pressure estimate in $L^{\infty}(\bL^2)$.
For the uniform estimates, we sketch a proof below. We note that in \cite{S95}, the error has been split into two:
$$\bu-\bue= (\bu-\bv)+(\bv-\bue)=\bxi+\bta,$$
where, $\bv$ is the solution of the linear penalized problem \cite[(4.1)-(4.2)]{S95}:
\begin{align*}
\bv_t-\nu\Delta\bv+\nabla\gamma = \f - B(\bu,\bu),\\
\nu\nabla\cdot \bv + \ve \gamma = 0, \quad \bv(0)=\bu_0,
\end{align*}
where $\bu$ is the solution of NSEs (\ref{om04})-(\ref{ibc04}) and $B(\bu,\bu)=(\bu\cdot\nabla)\bu$.
And the estimates of $\bxi$ are uniformly in time, see \cite[Lemma 4.1]{S95}. However, the estimates of $\bta$ grow exponentially in time (see, \cite[Theorem 4.1]{S95}) due to the use of the Gronwall's lemma. This can be avoided under the assumption of uniqueness condition (\ref{unique}).
We first note down the equation in $\bta$, see \cite[(4.7)]{S95}. 
\begin{equation}\label{s95bta}
\bta_t + \nu A_{\ve}\bta + \tilde{B}(\bue, \bxi+\bta)+ \tilde{B}(\bxi+\bta,\bu)=0.
\end{equation}
Take the inner product of (\ref{s95bta}) with $\bta$ to obtain 
\begin{equation}\label{s95bta1}
 \frac{1}{2}\frac{d}{dt}(\|\bta\|^2) + \nu \|A_{\ve}^{\frac{1}{2}}\bta\|^2 + \tb(\bue, \bxi+\bta,\bta)+ \tb(\bxi+\bta,\bu,\bta)=0.
\end{equation}
The nonlinear terms can be bounded by using the uniqueness condition (\ref{unique}) and the bounds \cite[(2.4)]{S95} as
\begin{align*}
\tb(\bue, \bxi+\bta,\bta)+ \tb(\bxi+\bta,\bu,\bta) &= \tb(\bue,\bxi,\bta)+\tb(\bxi,\bu,\bta)+\tb(\bta,\bu,\bta)\nonumber\\
& \le C(\|\Delta\bue\|+\|\Delta\bu\|)\|\bxi\|\|\nabla\bta\| + N \|\nabla\bu\|\|\nabla\bta\|^2 \nonumber\\
& \le C(\|\Delta\bue\|^2+\|\Delta\bu\|^2)\|\bxi\|^2 + \frac{\nu}{2}\|A_{\ve}^{\frac{1}{2}}\bta\|^2  + N  \|\nabla\bu\|\|\nabla\bta\|^2.
\end{align*}
Use the above estimate in (\ref{s95bta1}). Then, multiply by $e^{2\alpha t}$ and integrate with respect to time to arrive at
\begin{align}\label{s9uni}
  e^{2\alpha t}\|\bta(t)\|^2  + \nu\int_0^t e^{2\alpha s} \|A_{\ve}^{\frac{1}{2}}\bta(s)\|^2 ds  
  \le & \|\bta(0)\|^2 + 2\alpha \int_0^t e^{2\alpha s}\|\bta(s)\| ds  + N \int_0^t e^{2\alpha s}\|\nabla\bu(s)\|\|\nabla\bta(s)\|^2 ds \\
  &+ C \int_0^t e^{2\alpha s}(\|\Delta\bue(s)\|^2+\|\Delta\bu(s)\|^2)\|\bxi(s)\|^2ds.
\end{align}
The last term on the right hand side of (\ref{s9uni}) can be written as
\begin{align}\label{s9uni2}
\int_0^t e^{2\alpha s}(\|\Delta\bue(s)\|^2+\|\Delta\bu(s)\|^2)\|\bxi(s)\|^2ds 
&\le \|\bxi(t)\|_{L^\infty(\bL^2)}\int_0^t e^{2\alpha s}(\|\Delta\bue(s)\|^2+\|\Delta\bu(s)\|^2)ds \nonumber\\
& \le Ce^{2\alpha t} \|\bxi(t)\|_{L^\infty(\bL^2)}^2.
\end{align}
Now, we rewrite the last term on the left hand side of (\ref{s9uni}) using $\|A_{\ve}^{\frac{1}{2}}\bta\|^2= \|\nabla\bta\|^2+\frac{1}{\ve}\|\nabla\cdot\bta\|^2$. Then use (\ref{s9uni2}) and multiply the above estimate by $e^{-2\alpha t}$ to obtain
\begin{align*}
  \|\bta(t)\|^2  + e^{-2\alpha t}\int_0^t \big(\nu-2N \|\nabla\bu\|\big) e^{2\alpha s} \|\nabla\bta(s)\|^2 ds  + \frac{\nu}{\ve} e^{-2\alpha t}\int_0^t  e^{2\alpha s}\|\nabla\cdot\bta(s)\|^2 ds \\
  \le e^{-2\alpha t}\|\bta(0)\|^2 + 2\alpha e^{-2\alpha t} \int_0^t e^{2\alpha s}\|\bta(s)\|^2 ds + C  \|\bxi(t)\|_{L^\infty(\bL^2)}^2.
\end{align*}
Take limit on both sides as $t\to \infty$. Using 
$\overline{\lim}_{t\to\infty}\|\nabla\bu\| \le \frac{1}{\nu}\|\f\|_{L^\infty(0,\infty;\bH^{-1}(\Omega))}$,
it follows that
\begin{equation*}
\big(\nu-\frac{2N}{\nu}\|\f\|_{L^\infty(0,\infty;\bH^{-1}(\Omega))}\big)\overline{\lim_{t\to\infty}} \|\nabla \bta(t)\|^2  \le  C  \overline{\lim_{t\to\infty}} \|\bxi(t)\|_{L^\infty(\bL^2)}^2.
\end{equation*}
Since $ 2N\nu^{-2} \|\f\|_{L^\infty(0,\infty;\bH^{-1}(\Omega))}<1$, from (\ref{unique}), we conclude the following
\begin{equation*}
\overline{\lim_{t\to\infty}} \|\bta(t)\|  \le C \overline{\lim_{t\to\infty}} \|\nabla \bta(t)\|  \le  C\overline{\lim_{t\to\infty}} \|\bxi(t)\|_{L^\infty(\bL^2)}.
\end{equation*}
Combining the estimate of $\bxi$ from \cite[Lemma 4.1]{S95}, with the above estimate, we finally obtain
\begin{equation*}
\overline{\lim_{t\to\infty}} \|(\bu-\bue)(t)\|  \le  C\ve t^{-\frac{1}{2}}.
\end{equation*}
Using the above estimate, we can easily find the estimates of $\|\nabla(\bu-\bue)(t)\|$ and $\|(p-\pe)(t)\|$, which completes the rest of the proof.
\end{proof}

\section{Galerkin Finite Element Method}
\se

This section deals with the finite element Galerkin approximations to the penalized problem (\ref{pom}) or (\ref{pomu}) and some \textit{a priori} bounds for the semidiscrete problem.\\
Let $\mathcal{T}_h=\{K\}$ be a shape regular triangulation of the polygonal domain $\bar{\Omega}$ into closed subset $K$, triangles or quadrilaterals with $h = \max_{K}h_{K}$, where $h_{K}$ is the diameter of $K$. 

Let  $\bH_h$ and $L_h$ be two families of finite element spaces of $\bH_0^1 $ and $L^2/\R$, respectively, approximating the velocity vector and the pressure. It is assumed that the spaces $\bH_h$ and $L_h$ comprise of piecewise polynomial of degree at most $m$ and $m-1$ ($m\ge 2$), respectively. Assume that the following approximation properties are satisfied for the spaces $\bH_h$ and $L_h$: \\
%------------------------------------------------------------------------------(B1)
${\bf (B1)}$ For each $\bw \in\bH_0^1 \cap \bH^{m+1} $ and $ q \in H^m/\R$ with $m\ge 1$, there exist approximations $i_h w \in \bH_h $ and $ j_h q \in L_h $ such that
$$ \|\bw-i_h\bw\|+ h \| \nabla (\bw-i_h \bw)\| \le Ch^{j+1} \| \bw\|_{j+1}, ~~~~\| q - j_h q \| \le Ch^j \| q\|_{j},~~~ 0\le j \le m. $$
Further, we assume that inverse inequality holds for $\bw_h\in\bH_h$:
\begin{equation*} %\label{inv.hypo}
 \|\nabla \bw_h\| \leq  Ch^{-1} \|\bw_h\|.
\end{equation*}
We consider now the discrete analogue of the weak formulations (\ref{wfpom}) and (\ref{wfpp}): Find $(\bueh,\peh)$ in
$\bH_h\times L_h$ satisfying
\begin{equation}\left\{ \begin{array}{r}\label{dwfpom}
 (\bueht,\bphi_h)+\nu a(\bueh,\bphi_h)+\tb(\bueh,\bueh,\bphi_h) -(\peh,\nabla\cdot\bphi_h)=(\f,\bphi_h),~~\forall~\bphi_h\in\bH_h, \\
 \nu(\nabla\cdot\bueh,\chi_h)+\ve(\peh,\chi_h)=0,~~\forall~\chi_h\in L_h.
\end{array}\right.
\end{equation}
Choose $\chi_h=\nabla\cdot\bueh\in L_h$ in the second equation of (\ref{dwfpom}) and use this in the first equation to arrive at
\begin{equation}\label{dwfpomu}
 (\bueht,\bphi_h)+\nu a_{\ve}(\bueh,\bphi_h)+\tb(\bueh,\bueh,\bphi_h) =(\f,\bphi_h) ,~~\forall~\bphi_h\in\bH_h.
\end{equation}

\noindent For continuous dependence of the discrete pressure $\peh (t) \in L_h$ on the discrete velocity $\bueh(t) \in {\bf H}_h$, we assume the following discrete inf-sup (LBB) condition: \\
\noindent%----------------------------------------------------------------------(B2')
${\bf (B2^\prime)}$  For every $q_h \in L_h$, there exists a non-trivial function $\bphi_h \in \bH_h$ such that
$$ |(q_h, \nabla\cdot \bphi_h)| \ge C \|\nabla \bphi_h \|\| q_h\|, $$
where  the constant $C>0$ is independent of $h$. \\
Moreover, we also assume that the following approximation property holds true for ${\bf H}_h $. \\
\noindent%----------------------------------------------------------------------(B2)
${\bf (B2)}$ For every $\bw \in \bH_{0}^{1} \cap \bH^{m+1}, $ there exists an approximation $r_h \bw \in {\bH_h}$ such that
$$ \|\bw-r_h\bw\|+h \| \nabla (\bw - r_h \bw) \| \le Ch^{j+1} \|\bw\|_{j+1},~~~0\le j\le m. $$
Based on it, $L^2$-projection $P_h$,  defined as $P_h:\bL^2\to\bH_h$, can be derived, satisfying the following properties \cite{HR90}:
\begin{align*}   %\label{phpn01}
%\|\bphi-P_h \bphi\|+h\|\nabla P_h\bphi\| & \le Ch^{j}\|\bphi\|_{j},\quad\quad \quad\forall~\bphi\in \bH_h,     \\
    \|\bphi-P_h\bphi\|+h\|\nabla(\bphi-P_h\bphi)\| & \le Ch^{j+1}\|\bphi\|_{j+1},\quad\forall~
       \bphi\in\bH_0^1(\Omega)\cap\bH^{m+1}(\Omega). 
\end{align*}
Now we define the discrete operator $\Delta_h:\bH_h\to \bH_h$ through the bilinear form $a(\cdot,\cdot)$ as
\begin{equation*}
a(\bv_h,\bphi_h)=(-\Delta_h\bv_h,\bphi_h),~~~\forall~\bv_h,\bphi_h\in \bH_h.
\end{equation*}
Next we define the discrete analogue $A_{\ve h} : \bH_h\to\bH_h$ of $A_{\ve}$ satisfying 
\begin{align}\label{Aeh}
(A_{\ve h}\bv_h,\bphi_h) = a_\ve(\bv_h,\bphi_h) = a(\bv_h,\bphi_h)+\frac{1}{\ve}(\nabla\cdot\bv_h,\nabla\cdot\bphi_h),~\forall \bv_h,\bphi_h\in\bH_h.
\end{align}
Note that, $A_{\ve h}$ is a self-adjoint and  positive definite operator. With additional assumptions (\textbf{B1}) and (\textbf{B2}), the operator $A_{\ve h}$ mimics the estimates presented in the Lemmas \ref{Aep} and \ref{Aep1}.
\begin{lemma}\label{Aeph}%------------------------------------------A_{\ve} prop.
There exists a constant $c_0>0$ such that for $\ve>0$ sufficiently small, the following estimates hold:
\begin{align*}
& \|\Delta_h\bv_h\| \le  c_0\|A_{\ve h}\bv_h\|, \quad\forall~\bv_h\in \bH_h, \\
& \|\nabla\bv_h\| \le  c_0\|A_{\ve h}^{\frac{1}{2}}\bv_h\|, \quad\forall~\bv_h\in \bH_h, \\
& \|A_{\ve h}^{-\frac{r}{2}}\bv_h\| \le  c_0\|\bv_h\|_{-r},\quad\forall~\bv_h\in\bH_h,~~r\in \{1,2\}.
\end{align*}
\end{lemma}
\begin{proof}
Let $\bv_h\in\bH_h$ and $A_{\ve h}\bv_h=\g$. Take $q_h=-\frac{1}{\ve}\nabla\cdot\bv_h$, then (\ref{Aeh}) can be written as
\begin{align*}
(\nabla\bv_h,\nabla\bphi_h)- (q_h,\nabla\cdot\bphi_h)= (\g,\bphi_h),~~\forall \bphi_h\in\bH_h,\\
(\nabla\cdot\bv_h,\psi_h) + \ve(q_h,\psi_h)=0 ,~~\forall \psi_h\in L_h.
\end{align*}
From regularity estimate, one can find that (see, \cite[(1.20)]{BGT88})
\begin{align*}%\label{reg}
\|\Delta_h\bv_h\| + \|\nabla q_h\| \le c_0 \|\g\| + \ve c_0 \|\nabla q_h\|.
\end{align*}

\noindent
Now choose $\ve$ sufficiently small such that $c_0\ve<1$, then we conclude the first result. 
For the second result, we choose $\bphi_h=\bv_h$ in (\ref{Aeh}) and arrive at 
\begin{align*}
\|A_{\ve h}^{\frac{1}{2}}\bv_h\|^2 = \|\nabla\bv_h\|^2 + \frac{1}{\ve}\|\nabla\cdot\bv_h\|^2 \ge \|\nabla\bv_h\|^2 .
\end{align*}
For the third one, let $\bw_h$ be the solution of $A_{\ve h}^{\frac{r}{2}}\bw_h=A_{\ve h}^{-\frac{r}{2}}\bv_h$.
\begin{align*}
\|A_{\ve h}^{-\frac{r}{2}}\bv_h\|^2 = (A_{\ve h}^{\frac{r}{2}}\bw_h, A_{\ve h}^{-\frac{r}{2}}\bv_h) = (\bw_h,\bv_h) \le \|\bw_h\|_{r}\|\bv_h\|_{-r}
 \le c_0\|A_{\ve h}^{\frac{r}{2}}\bw_h\|\|\bv_h\|_{-r} \le  c_0\|A_{\ve h}^{-\frac{r}{2}}\bv_h\|\|\bv_h\|_{-r}.
\end{align*}
Cancelling one $\|A_{\ve h}^{-\frac{r}{2}}\bv_h\|$ from both sides completes the rest of the proof.
\end{proof}

\noindent
Before we proceed further, we look at some standard estimates of the nonlinear term $b(\cdot,\cdot,\cdot)$ (see \cite{HR90, GD15}), which use H\"older's inequality and the following discrete Ladyzhenskaya's inequality \cite{T84}, for all $\bphi_h\in\bH_h$:
\begin{align*}%\label{lady}
& \|\bphi_h\|_{L^4} \le C \|\bphi_h\|^{\frac{1}{2}}\|\nabla\bphi_h\|^{\frac{1}{2}},\\
& \|\nabla\bphi_h\|_{L^4} \le C \|\nabla\bphi_h\|^{\frac{1}{2}}\|\Delta_h\bphi_h\|^{\frac{1}{2}},\\
& \|\bphi_h\|_{L^\infty} \le C \|\bphi_h\|^{\frac{1}{2}}\|\Delta_h\bphi_h\|^{\frac{1}{2}}.
\end{align*} 
\begin{lemma} \label{trilinear}
Suppose the conditions $({\bf A1}), ({\bf B1})$ and $({\bf B2})$ are satisfied. Then there exists a constant $C>0$ such that for all $\bv_h, \bw_h, \bphi_h \in \bH_h$, the trilinear form $b(\cdot,\cdot,\cdot)$ satisfies the following properties:
\[
|\tb(\bv_h,\bw_h,\bphi_h)|= C
\begin{cases}
\|\bv_h\|^{\frac{1}{2}} \|\Delta_h\bv_h\|^{\frac{1}{2}} \|\nabla\bw_h\| \|\bphi_h\|, \\
\|\bv_h\|^{\frac{1}{2}}\|\nabla\bv_h\|^{\frac{1}{2}} \|\nabla\bw_h\|^{\frac{1}{2}} \|\Delta_h\bw_h\|^{\frac{1}{2}} \|\bphi_h\|, \\
\|\bv_h\| \|\nabla\bw_h\| \|\bphi_h\|^{\frac{1}{2}}\|\Delta_h\bphi_h\|^{\frac{1}{2}}, \\
\|\bv_h\| \|\nabla\bw_h\|^{\frac{1}{2}}\|\Delta_h\bw_h\|^{\frac{1}{2}} \|\bphi_h\|^{\frac{1}{2}}\|\nabla\bphi_h\|^{\frac{1}{2}},\\
\|\bv_h\|^{\frac{1}{2}}\|\nabla\bv_h\|^{\frac{1}{2}} \|\nabla\bw_h\| \|\bphi_h\|^{\frac{1}{2}}\|\nabla\bphi_h\|^{\frac{1}{2}},
\end{cases}
\]
\end{lemma}

%\subsection{A priori estimates}
\noindent 
We next look at the {\it a priori} estimates of the discrete penalized solution $\bueh$. Similar to the continuous case (like in Lemma \ref{pap}), these estimates can be easily verified.
\begin{lemma}\label{paph}
 Apart from the assumptions of the Lemma \ref{pap}, we assume that $({\bf B1})$ and $({\bf B2})$ hold. Then, there exists constant $C$ independent of $\ve$ and $h$ such that for $\bueh(0)=P_h\bu_{\ve 0}$
\begin{eqnarray*}
 \|\bueh(t)\|^2+e^{-2\alpha t}\int_0^t e^{2\alpha s} \|A_{\ve h}^{\frac{1}{2}}\bueh(s)\|^2 ds   &\le& C,  \\
 \tau^{r}(t)\|A_{\ve h}^{\frac{r+1}{2}}\bueh(t)\|^2 + e^{-2\alpha t}\int_0^t e^{2\alpha s}\|A_{\ve h} \bueh(s)\|^2 ds &\le & C,~~~ r \in \{0,1\}, \\
 e^{-2\alpha t}\int_0^t \tau^{r}(s) e^{2\alpha s}\|A_{\ve h}^{\frac{r-2}{2}}\buehss(s)\|^2 ds &\le & C,~~~ r \in\{0,1,2\},
\end{eqnarray*}
hold, where $\tau(t)=\min\{t,1\}$.
\end{lemma}
\begin{proof}
Choose $\bphi_h=\bueh$ in (\ref{dwfpomu}) and use the Cauchy-Schwarz inequality and the Poincar\'e inequality with Lemma \ref{Aeph} ($\|\bueh\|^2 \le \frac{1}{\lambda_1}\|\nabla\bueh\|^2\le \frac{c_0^2}{\lambda_1}\|A_{\ve h}^{\frac{1}{2}}\bueh\|^2$) to find that
\begin{align}\label{pap000}
 \frac{d}{dt}\|\bueh\|^2  + \nu \|A_{\ve h}^{\frac{1}{2}}\bueh\|^2  \le  \frac{c_0^2}{\nu\lambda_1}\|\f\|^2.
\end{align}
Note that the non-linear term vanishes due to (\ref{tbp1}). Now multiply by $e^{2\alpha t}$ and integrate from $0$ to $t$ to obtain
\begin{align}\label{pap001}
 e^{2\alpha t}\|\bueh(t)\|^2+(\nu-\frac{2c_0^2\alpha}{\lambda_1})\int_0^t e^{2\alpha s}\|A_{\ve h}^{\frac{1}{2}}\bueh(s)\|^2ds \le 
\|P_h\bu_{\ve 0}\|^2+\frac{c_0^2(e^{2\alpha t}-1)}{2\alpha\nu\lambda_1}\|\f\|^2_{\infty}.
\end{align}
With $0< \alpha < \frac{\nu\lambda_1}{2c_0^2}$, we have $(\nu-\frac{2c_0^2\alpha}{\lambda_1})>0$. Multiply through out by $e^{-2\alpha t}$ to conclude the first proof. 
Now, we integrate (\ref{pap000}) with respect to time from $t$ to $t+T$ for any $T>0$, we have
\begin{align}\label{pap0015}
\|\bueh(t+T)\|^2+\nu\int_t^{t+T} \|A_{\ve h}^{\frac{1}{2}}\bueh(s)\|^2ds \le 
\|\bueh(t)\|^2 + \frac{c_0^2T}{\nu\lambda_1}\|\f\|^2_{\infty}.
\end{align}
For the second estimate, choose $\bphi_h=A_{\ve h}^{r+1}\bueh$, $r \in \{0,1\}$ in (\ref{dwfpomu}). When $r=0$, we find that
\begin{align}\label{pap01}
 \frac{1}{2}\frac{d}{dt} \|A_{\ve h}^{\frac{1}{2}}\bueh\|^2 +\nu\|A_{\ve h}\bueh\|^2 = (\f,A_{\ve h}\bueh)  -\tb(\bueh,\bueh,A_{\ve h}\bueh).
\end{align}
We use Lemma \ref{trilinear} with Lemma \ref{Aeph}, the Young's inequalities to bound the nonlinear term as
\begin{align} \label{pap011}
 \tb(\bueh,\bueh,A_{\ve h}\bueh) 
%&\le \|\bueh\|^{\frac{1}{2}}\|\nabla\bueh\| \|\Delta_h\bueh\|^{\frac{1}{2}} \|A_{\ve h}\bueh\|  \\
\le C  \|\bueh\|^{\frac{1}{2}}\|A_{\ve h}^{\frac{1}{2}}\bueh\|\|A_{\ve h}\bueh\|^{3/2} 
\le C \|\bueh\|^2\|A_{\ve h}^{\frac{1}{2}}\bueh\|^4 + \frac{\nu}{4}\|A_{\ve h}\bueh\|^2. 
\end{align}
Substitute the above estimate in (\ref{pap01}) to find that
\begin{align}\label{pap02}
 \frac{d}{dt} (\|A_{\ve h}^{\frac{1}{2}}\bueh\|^2) +\nu \|A_{\ve h}\bueh\|^2 
 \le  C  \big(\|\bueh\|^2\|A_{\ve h}^{\frac{1}{2}}\bueh\|^2\big) \|A_{\ve h}^{\frac{1}{2}}\bueh\|^2 + \frac{2}{\nu}\|\f\|^2).
\end{align}
We now apply the uniform Gronwall's Lemma (Lemma \ref{ugl}) in (\ref{pap02})  and use (\ref{pap001}) and (\ref{pap0015}) to conclude that $\|A_{\ve h}^{\frac{1}{2}}\bueh(t+T)\|^2$ is uniformly bounded with respect to $t$ for all $t\ge 0$, that is,  $\|A_{\ve h}^{\frac{1}{2}}\bueh(t)\|^2$ is uniformly bounded on $[T,\infty)$. Precisely
\begin{equation}\label{pap4}
  \|A_{\ve h}^{\frac{1}{2}}\bueh(t)\|^2  \le C, \quad \forall t\ge T.
\end{equation}
For $0\le t\le T$, we use the classical Gronwall's lemma \cite{S95} in (\ref{pap02}) and obtain 
\begin{equation}\label{pap41}
  \|A_{\ve h}^{\frac{1}{2}}\bueh(t)\|^2  \le C, \quad \text{for}~~ 0\le t\le T.
\end{equation}
Finally,  multiply (\ref{pap02}) by $e^{2\alpha t}$ and integrate with respect to time from $0$ to $t$ and use the estimates (\ref{pap001}), (\ref{pap4}) and (\ref{pap41}) to complete the second proof when $r=0$. 
For $r=1$, we need some intermediate estimate. First we take $\bphi_h=e^{2\alpha t}\bueht$ with $\hbueh=e^{\alpha t}\bueh$ in (\ref{dwfpomu}) to obtain
\begin{align}\label{ppap01}
 \frac{\nu}{2}\frac{d}{dt} \|A_{\ve h}^{\frac{1}{2}}\hbueh\|^2+  \|\hbueht\|^2= \alpha \nu\|A_{\ve h}^{\frac{1}{2}}\hbueh\|^2 
  + (\hf,\hbueht) -e^{2\alpha t}\tb(\bueh,\bueh,\bueht).
\end{align}
We can estimate the nonlinear term on the right hand side of (\ref{ppap01}) using Lemma \ref{trilinear} and 
integrate both sides with respect to time to find that
\begin{align*}%\label{ppap02}
 \nu \|A_{\ve h}^{\frac{1}{2}}\hbueh\|^2 + \int_0^t\|\hbueht\|^2 ds
  \le C \bigg[\int_0^t \Big(\|A_{\ve h}^{\frac{1}{2}}\hbueh\|^2 + \|\hf\|^2+ \|A_{\ve h}^{\frac{1}{2}}\bueh\|^2\|A_{\ve h}\hbueh\|^2 \big) ds\bigg].
\end{align*}
Now a use of (\ref{pap001}) and (\ref{pap4}) lead us to the intermediate estimate. 
\begin{align} \label{ppap021}
  \|A_{\ve h}^{\frac{1}{2}}\bueh(t)\|^2 + e^{-2\alpha t}\int_0^t e^{2\alpha s}\|\buehs(s)\|^2 ds \le C.
\end{align}
We now differentiate (\ref{dwfpomu}) with respect to time and deduce that
\begin{align}\label{dwfpomut}
 (\buehtt,\bphi_h)+\nu a_{\ve}(\bueht,\bphi_h) = (\f_t,\bphi_h) -\tb(\bueht,\bueh,\bphi_h)-\tb(\bueh,\bueht,\bphi_h), \quad \forall \bphi_h\in \bH_h.
\end{align}
Take $\bphi_h=\sigma(t)\buehtt$, where $\sigma(t)=\tau(t) e^{2\alpha t}$ and use Lemma \ref{trilinear} with Lemma \ref{Aeph}, the Cauchy-Schwarz inequality to reach at 
\begin{align*}
\frac{d}{dt}(\sigma(t)\|\bueht\|^2) + \nu \sigma(t)\|A_{\ve h}^{\frac{1}{2}}\bueht\|^2 \le C e^{2\alpha t}\|\bueht\|^2 +C\sigma(t)\Big(\|\f_t\|^2+ \|\bueht\|^2\|A_{\ve h}^{\frac{1}{2}}\bueh\|^2\Big).
\end{align*}
Integrate  with respect to time and use (\ref{ppap021}), (\ref{pap4}) and (\ref{pap41}) to obtain
\begin{align} \label{ppap022}
  \tau(t)\|\bueht(t)\|^2 + \nu  e^{-2\alpha t}\int_0^t \sigma(s)\|A_{\ve h}^{\frac{1}{2}}\buehs(s)\|^2 ds \le C.
\end{align}
Now we are in position to complete the proof of the second estimate when $r=1$. Set $\bphi_h=A_{\ve h}\bueh$ in (\ref{dwfpomu}) and rewrite it and use (\ref{pap011}) and the Cauchy-Schwarz inequality to arrive at
\begin{align*}
\nu \|A_{\ve h}\bueh\|^2 &= (\f,A_{\ve h}\bueh)-(\bueht,A_{\ve h}\bueh)-\tb(\bueh,\bueh,A_{\ve h}\bueh)\\
 &\le  C  \big( \|\f\|^2 +\|\bueht\|^2  + \|\bueh\|^2\|A_{\ve h}^{\frac{1}{2}}\bueh\|^4\big) + \frac{\nu}{2}\|A_{\ve h}\bueh\|^2.
\end{align*}
Multiply by $\tau(t)$ and use (\ref{pap001}), (\ref{pap4}), (\ref{pap41}) and (\ref{ppap022}) to complete the second proof.
\\
For the third estimate, choose $\bphi_h=e^{2\alpha t}A_{\ve h}^{-2}\buehtt$ in (\ref{dwfpomut}). Then, we use (\ref{tbp1}) and Lemma \ref{trilinear} with Lemma \ref{Aeph}  and the Cauchy-Schwarz inequality to bound the terms on right hand side as
\begin{align*}
 |(\f_t,A_{\ve h}^{-2}\buehtt) & - \tb(\bueht,\bueh,A_{\ve h}^{-2}\buehtt)-\tb(\bueh,\bueht,A_{\ve h}^{-2}\buehtt)| \\
&\le |(\f_t,A_{\ve h}^{-2}\buehtt)|+ |\tb(\bueht,\bueh,A_{\ve h}^{-2}\buehtt)|+ |\tb(\bueh,A_{\ve h}^{-2}\buehtt,\bueht)| \\
&\le C (\|\f_t\|^2 + \|\bueht\|^2\|A_{\ve h}^{\frac{1}{2}}\bueh\|^2) + \frac{1}{2}\|A_{\ve h}^{-1}\buehtt\|^2.
\end{align*}
Then we arrive at  
\begin{align*}%\label{dwfpomut1}
 e^{2\alpha t}\|A_{\ve h}^{-1}\buehtt\|^2 + \frac{\nu}{2}\frac{d}{dt}(e^{2\alpha t}\|A_{\ve h}^{-\frac{1}{2}}\bueht\|^2) \le \nu\alpha e^{2\alpha t}\|A_{\ve h}^{-\frac{1}{2}}\bueht\|^2 + C e^{2\alpha t}\Big( \|\bueht\|^2\|A_{\ve h}^{\frac{1}{2}}\bueh\|^2+\|\f_t\|^2\Big).
\end{align*}
Integrate both sides with respect to time from $0$ to $t$. Now, a use of the estimates obtained above results in the case $r=0$. For  $r=1$ and $r=2$, we take $\bphi_h=\sigma (t)A_{\ve h}^{-1}\buehtt$ and $\bphi_h=\sigma^{2}(t)\buehtt$ with $\sigma^{2}(t)=(\tau)^2(t) e^{2\alpha t}$, respectively, and do similar analysis as above. This completes the rest of the proof.
\end{proof}

\section{Error Estimates for the Semidiscrete Problem}
\se

This section deals with the error analysis of the finite element Galerkin approximation for the penalized system (\ref{wfpp}). 
Here, the goal is to provide optimal $L^{\infty}(L^2)$ error estimates for the velocity and the pressure when $\bu_{\ve 0}\in \bH_0^1$. The main result of this section is as follows.
%------------------------------------------------PENALIZED ERROR ESTIMATE
\begin{theorem}\label{perrest}
Let the conditions $({\bf A1})$,$({\bf A2})$, $({\bf B1})$ and $({\bf B2})$ be satisfied. Further, let the discrete initial velocity $\bueh(0)=P_h\bu_{\ve 0},$ where $\bu_{\ve 0}\in \bH_0^1(\Omega).$ 
Then, there exists a positive constant $C$ such that for $t>0$
\begin{equation*}%\label{perrest0}
 \|(\bue-\bueh)(t)\|+h(\|\nabla(\bue-\bueh)(t)\|+\|(\pe-\peh)(t)\|) \le K(t) h^{m+1}t^{-\frac{m}{2}}.
\end{equation*}
where, $K(t)=Ce^{Ct}$. 
Under the additional condition (\ref{unique}), the estimates are uniformly in time, that is $K(t)=C$.
\end{theorem}

\noindent
The proof of the theorem is realized via a numbers of lemmas. However before we visit the lemmas, we need to develop some preliminary tools, with which we begin our discussion. \\ 
Let $\f\in \bL^2$. We define the linear inverse operators $A_{\ve}^{-1}:\bL^2\to\bH_0^1$ and $A_{\ve h}^{-1}:\bH_h\to\bH_h$ satisfying
\begin{align*}
& a_\ve(A_{\ve}^{-1}\f,\bphi) = (\nabla A_{\ve}^{-1}\f,\nabla\bphi)+\frac{1}{\ve}(\nabla \cdot A_{\ve}^{-1}\f,\nabla\cdot\bphi)=(\f,\bphi), \quad\forall \bphi\in\bH_0^1  \\ %\label{AeI}
& a_\ve(A_{\ve h}^{-1}P_h\f,\bphi_h) = (\nabla A_{\ve h}^{-1}P_h\f,\nabla\bphi_h)+\frac{1}{\ve}(\nabla \cdot A_{\ve h}^{-1}P_h\f,\nabla\cdot\bphi_h)=(P_h\f,\bphi_h), \quad\forall \bphi_h\in\bH_h.  %\label{AehI}
\end{align*} 
Following the work of Heywood and Rannacher \cite{HR82}, we have the following result.
\begin{proposition}\label{pro1}
The map $A_{\ve h}^{-1}P_h A_{\ve}:\bH_0^1\cap\bH^2\to \bH_h$ satisfies the following estimate:
\begin{align*}
\|\bv- A_{\ve h}^{-1}P_h A_{\ve}\bv\| + h\|\nabla(\bv- A_{\ve h}^{-1}P_h A_{\ve}\bv)\| \le Ch^{m+1}\|A_{\ve}^{\frac{m+1}{2}}\bv\|.
\end{align*}
\end{proposition}

\noindent
We next consider the penalized steady Stokes problem: For a given function $\g\in\bL^2$, let $\bv\in\bH_0^1\cap\bH^2$, $q\in H^1/R$ be the unique solution of 
\begin{align*}
a(\bv,\bphi)-(q,\nabla\cdot\bphi) = (\g,\bphi),~~\forall \bphi\in\bH_0^1,\\
(\nabla\cdot\bv, \chi) + \ve(q,\chi)=0,~~\forall \chi\in L^2.
\end{align*}
Getting rid of $q$, we obtain
\begin{align}\label{stf}
a(\bv,\bphi)+ \frac{1}{\ve}(\nabla\cdot\bv,\nabla\cdot\bphi) = (\g,\bphi),~~\forall \bphi\in\bH_0^1.
\end{align}
The finite element approximation $\bv_h\in\bH_h$ of $\bv$ satisfies the following equation
\begin{align}\label{dstf}
a(\bv_h,\bphi_h)+ \frac{1}{\ve}(\nabla\cdot\bv_h,\nabla\cdot\bphi_h) = (\g,\bphi_h),~~\forall \bphi_h\in\bH_h.
\end{align}
\begin{lemma}\label{psterr} %penalized stokes error
Let the conditions $({\bf A1})$,$({\bf A2})$, 
$({\bf B1})$ and $({\bf B2})$ be satisfied. Then, there exists a positive constant $C$ such that 
\begin{equation*}
 \|\bv-\bv_h\|+h\|\nabla(\bv-\bv_h)\| \le C h^{m+1}\|A_{\ve}^{\frac{m+1}{2}}\bv\|.
\end{equation*}
\end{lemma}
\begin{proof}
The main idea of the proof is adapted from the paper of Heywood and Rannacher \cite{HR82}. From (\ref{stf}) and (\ref{dstf}), we have the following error equation 
\begin{align}\label{sterr1}
a(\bv-\bv_h,\bphi_h)+ \frac{1}{\ve}(\nabla\cdot(\bv-\bv_h),\nabla\cdot\bphi_h) = 0,~~\forall \bphi_h\in\bH_h.
\end{align}
We choose $\bphi_h=i_h\bv-\bv_h\in\bH_h$ in the above equation to obtain
\begin{align*}
\|\nabla(\bv-\bv_h)\|^2+ \frac{1}{\ve}\|\nabla\cdot(\bv-\bv_h)\|^2 
= a(\bv-\bv_h,\bv-i_h\bv)+ \frac{1}{\ve}(\nabla\cdot(\bv-\bv_h),\nabla\cdot(\bv-i_h\bv)).
\end{align*}
A use of the Cauchy-Schwarz inequality with (\textbf{B1}) and Lemma \ref{Aep} yields
\begin{align}\label{sterr11}
\|\nabla(\bv-\bv_h)\|^2+ \frac{1}{\ve}\|\nabla\cdot(\bv-\bv_h)\|^2 
&\le C (\|\nabla(\bv-i_h\bv)\|^2 + \frac{1}{\ve}\|\nabla\cdot(\bv-i_h\bv)\|^2) \nonumber\\
&\le Ch^{2m}(\|\bv\|_{m+1}^2+\frac{1}{\ve}\|\nabla\cdot\bv\|_{m}^2) \nonumber\\
&\le Ch^{2m}\|A_{\ve}^{\frac{m+1}{2}}\bv\|^2.
\end{align}
To obtain the $L^2$ estimate, we need to consider a duality problem
\begin{align}\label{stdual}
A_{\ve}\bw=\bv-\bv_h, ~~ \bw|_\Omega=0,
\end{align}
satisfying
\begin{align}\label{stdualap}
\|A_{\ve}\bw\| \le C \|\bv-\bv_h\|.
\end{align}
We now take inner product on both sides of (\ref{stdual}) by $\bv-\bv_h$ and choose $\bphi_h=i_h\bw$ in (\ref{sterr1}) to obtain
\begin{align*}
\|\bv-\bv_h\|^2 & = a(\bv-\bv_h,\bw-i_h\bw)+ \frac{1}{\ve}(\nabla\cdot(\bv-\bv_h),\nabla\cdot(\bw-i_h\bw)).
\end{align*}
A use of the Cauchy-Schwarz inequality with the Young's inequality with $\delta>0$, approximation property (\textbf{B1}), Lemma \ref{Aeh} and (\ref{stdualap}) shows
\begin{align*}
\|\bv-\bv_h\|^2 
&\le Ch^{2m+2}\|A_{\ve}^{\frac{m+1}{2}}\bv\|^2 + \delta(\|\bw\|_2^2+\frac{1}{\ve}\|\nabla\cdot\bw\|_1^2)\\
&\le Ch^{2m+2}\|A_{\ve}^{\frac{m+1}{2}}\bv\|^2 + \delta\|A_{\ve}\bw\|^2 \\
&\le Ch^{2m+2}\|A_{\ve}^{\frac{m+1}{2}}\bv\|^2 + C\delta \|\bv-\bv_h\|^2.
\end{align*}
Now, a choice of $C\delta=\frac{1}{2}$ completes the rest of the proof.
\end{proof}

\subsection{Error Estimates for the Velocity}

\noindent 
In order to obtain the semidiscrete error estimates for the velocity, we denote $\eve=\bue-\bueh$ and subtract (\ref{dwfpomu}) from (\ref{wfpp}) to find the equation of the semi-discrete error $\eve$:
%
%-------------------------------PENALIZED PROBLEM: ERROR EQUATION; SEMI-DISCRETE
\begin{equation}\label{ppeesd}
 (\evet,\bphi_h)+\nu a_{\ve}(\eve,\bphi_h) =\tb(\bueh,\bueh,\bphi_h)-\tb(\bue,\bue,\bphi_h),\quad \forall\bphi_h\in\bH_h.
\end{equation}
By introducing an intermediate solution $\btveh$ which is a finite element Galerkin approximation to a linearized penalized NSEs, that is, $\btveh$ satisfies
%------------------------------------------------linear penalized NSE
\begin{equation}\label{lpo}
 (\btveht,\bphi_h)+\nu a_{\ve}(\btveh,\bphi_h) =(\f,\bphi_h)-\tb(\bue,\bue,\bphi_h)~~~~\forall \bphi_h\in\bH_h;
\end{equation}
we split $\eve$ as
$$ \eve := \bue - \bueh = (\bue - \btveh) + ( \btveh - \bueh)= \bxi_h + \bta_h. $$
Note that $\bxi_h$ is  the error committed by  approximating a linearized penalized NSEs and $\bta_h$ represents  the error due to the presence of non-linearity in the equation.  Below, we derive some estimates of $\bxi_h$. Subtracting (\ref{lpo}) from
(\ref{wfpp}), the equation in $\bxi_h$ is written as
\begin{equation}\label{pbxi}
 (\bxi_{ht},\bphi_h)+\nu a_{\ve}(\bxi_h,\bphi_h)=0,~~  \bphi_h\in\bH_h.
\end{equation}
%%----------------------------------------------L^2L^2 est of \bxi_h 
\begin{lemma}\label{l2epbxi}
Let the assumptions of Lemma \ref{paph} hold and $\btveh(t)\in\bH_h$ be a solution of (\ref{lpo}) with initial condition $\btveh(0)=P_h\bu_{\ve 0}$ and $\bue$ be a weak solution of (\ref{pomu}) with initial condition $\bu_{\ve 0} \in \bH_0^1$. Then, $\bxi_h$ satisfies
\begin{equation*}
 e^{-2\alpha t}\int_0^t \sigma^{m-1}(s)\|\bxi_h(s)\|^2 ~ds \le C h^{2m+2},\quad m\ge 1,
\end{equation*}
where, $\sigma^m(t)=\tau^m(t)e^{2\alpha t}$ and $\tau(t)=\min\{1,t\}$.
\end{lemma} 
\begin{proof}
We rewrite the equation (\ref{pomu}) and (\ref{lpo}) as
\begin{align*}
 A_{\ve}^{-1}\buet +\nu  \bue = A_{\ve}^{-1} (\f -\tilde B(\bue,\bue)),
\end{align*}
and
\begin{align*}
 A_{\ve h}^{-1}\btveht +\nu  \btveh =A_{\ve h}^{-1}P_h(\f -\tilde B(\bue,\bue)).
\end{align*}
In the view of above two equations, along with (\ref{pomu}) we have
\begin{align}\label{aveest0d}
A_{\ve h}^{-1}P_h\bxi_{ht} + \nu \bxi_h 
&= A_{\ve h}^{-1}(P_h\buet-\btveht) + \nu (\bue-\btveh) \nonumber\\
&= (A_{\ve h}^{-1}P_h-A_{\ve}^{-1})\buet + (A_{\ve}^{-1}\buet + \nu \bue)  - (A_{\ve h}^{-1}\btveht + \nu \btveh) \nonumber\\
&= (A_{\ve h}^{-1}P_h-A_{\ve}^{-1})\buet + A_{\ve}^{-1}(\f-\tilde B(\bue,\bue))  - A_{\ve h}^{-1}P_h(\f-\tilde B(\bue,\bue)) \nonumber\\
&= (A_{\ve h}^{-1}P_h-A_{\ve}^{-1})\buet + (A_{\ve}^{-1} - A_{\ve h}^{-1}P_h)(\buet - \nu A_{\ve}\bue) \nonumber\\
&=  - \nu (\bue - A_{\ve h}^{-1}P_h A_{\ve}\bue).
\end{align}
Taking inner product with $\bxi_h$ in above equation to arrive at
\begin{align}\label{aveest1d}
\frac{1}{2}\frac{d}{dt}\|A_{\ve h}^{-\frac{1}{2}}P_h\bxi_h\|^2 + \nu \|\bxi_h\|^2 
=  - \nu ((\bue - A_{\ve h}^{-1}P_h A_{\ve}\bue),\bxi_h).
\end{align}
A use of the Cauchy-Schwarz inequality, Proposition \ref{pro1}  in (\ref{aveest1d}) gives
\begin{align}\label{aveest3d}
 \frac{d}{dt}\|A_{\ve h}^{-\frac{1}{2}}P_h\bxi_h\|^2 + \nu \|\bxi_h\|^2 
\le   C h^{2m+2}\|A_{\ve}^{\frac{m+1}{2}}\bue\|^2.
\end{align}
We now multiply both sides of (\ref{aveest3d}) by $\sigma^{m-1}(t)$ and use the fact $\frac{d}{dt}\sigma^{m-1}(t)\le (m-1)\sigma^{m-2}(t)+ 2\alpha \sigma^{m-1}(t)$ to obtain
\begin{align*}
  \frac{d}{dt}(\sigma^{m-1}\|A_{\ve h}^{-\frac{1}{2}}P_h\bxi_h\|^2) - (m-1)\sigma^{m-2} \|A_{\ve h}^{-\frac{1}{2}}P_h\bxi_h\|^2 - 2\alpha \sigma^{m-1} \|A_{\ve h}^{-\frac{1}{2}}P_h\bxi_h\|^2 + \nu \sigma^{m-1}\|\bxi_h\|^2 \\
 \le   C h^{2m+2}\sigma^{m-1}\|A_{\ve}^{\frac{m+1}{2}}\bue\|^2.
\end{align*}
The third term on the left hand side can be combined with the forth term using the fact $ \|A_{\ve h}^{-\frac{1}{2}}P_h\bxi_h\|^2 \le c_0^2\|P_h\bxi_h\|_{-1}^2 \le \frac{c_0^2}{\lambda_1}\|\bxi_h\|^2$ and we finally 
integrate both sides from $0$ to $t$ and use $\|A_{\ve h}^{-\frac{1}{2}}P_h\bxi_h(0)\| = 0$, we deduce that
\begin{align}\label{aveest6d}
&\sigma^{m-1}(t)\|A_{\ve h}^{-\frac{1}{2}}P_h\bxi_h(t)\|^2 + \Big(\nu-\frac{2c_0^2\alpha}{\lambda_1}\Big) \int_0^t \sigma^{m-1}(s)\|\bxi_h(s)\|^2 ds \nonumber\\
&\qquad\le C h^{2m+2} \int_0^t \sigma^{m-1}(s)\|A_{\ve}^{\frac{m+1}{2}}\bue(s)\|^2 ds + (m-1)\int_0^t \sigma^{m-2}(s) \|A_{\ve h}^{-\frac{1}{2}}P_h\bxi_h(s)\|^2 ds.
\end{align}
With $0< \alpha < \frac{\nu\lambda_1}{2c_0^2}$, we have $(\nu-\frac{2c_0^2\alpha}{\lambda_1})>0$.
Now we prove the assertion of Lemma \ref{l2epbxi} in sequence of steps. First we take $m=1$, then the last term on the right hand side of (\ref{aveest6d}) vanishes and we obtain
\begin{equation*} %\label{aveestm1d}
 \int_0^t e^{2\alpha s}\|\bxi_h(s)\|^2 ~ds \le C h^{4} \int_0^t e^{2\alpha s}\|A_{\ve}\bue(s)\|^2 ds.
\end{equation*} 
For $m=2$, we have from (\ref{aveest6d})
\begin{equation}\label{aveestm2d}
 \int_0^t \sigma(s)\|\bxi_h(s)\|^2 ~ds \le C h^{6} \int_0^t \sigma(s)\|A_{\ve}^{3/2}\bue(s)\|^2 ds +  \int_0^t e^{2\alpha s}\|A_{\ve h}^{-\frac{1}{2}}P_h\bxi_h(s)\|^2.
\end{equation}
To find the estimate for the last term on the right hand side of (\ref{aveestm2d}), we take inner product with $A_{\ve h}^{-1}P_h\bxi_h$ in (\ref{aveest0d}) and argue as same as (\ref{aveest1d})-(\ref{aveest6d}), we find that
\begin{equation*} %\label{aveestm21d}
 \int_0^t e^{2\alpha s}\|A_{\ve h}^{-\frac{1}{2}}P_h\bxi_h(s)\|^2 ~ds \le C h^{6} \int_0^t e^{2\alpha s}\|A_{\ve}\bue(s)\|^2 ds.
\end{equation*}
For $m>2$, we can follow the same technique, which completes the rest of the proof.
\end{proof}

\noindent 
For optimal estimate of $\bxi_h$ in $L^{\infty}(\bL^2),$ we consider a projection $\Seh:\bH_0^1\to\bH_h$ satisfy
\begin{equation}\label{psve}
 a_{\ve}(\bue-\Seh\bue,\bphi_h) =0,~~~~\forall~\bphi_h\in\bH_h,
\end{equation}
for some fixed $\ve>0.$ We note that the above system, similar to \cite[(4.52)]{HR82}, has a positive definite operator $A_{\ve h}.$ Therefore, we can establish the well-posedness of the system (\ref{psve}) similar to \cite{HR82}.
 
\noindent
We now write
$$ \bxi_h=(\bue-\Seh\bue)+(\Seh\bue-\btveh) =:\bzeta+\btheta. $$
We are interested in the estimates of $\|\bzeta\|,\|\nabla\bzeta\|,$ as this is the first step towards obtaining the optimal estimate of $\bxi_h$. We present the following Lemma.
%%-----------------------------------------Penalized Stokes-Volterra estimates: \bue-\Seh\bue
\begin{lemma}\label{psvest0}
Suppose the assumptions of Lemma \ref{paph} are satisfied. Then, there exists a positive constant $C$ such that
\begin{equation*}%\label{psv01}
 \|\bzeta(t)\|+ h\|\nabla\bzeta(t)\| \le Ch^{m+1} \|A_{\ve}^{\frac{m+1}{2}}\bue(t)\|.
\end{equation*}
Moreover, the following estimate holds:
\begin{equation*}%\label{psv02}
 \|\bzeta_t(t)\|+ h\|\nabla\bzeta_t(t)\|\le Ch^{m+1}\|A_{\ve}^{\frac{m+1}{2}}\buet(t)\|.
\end{equation*}
\end{lemma}
\begin{proof}
Note that, $\bzeta=\bue-\Seh\bue,$ then the first result directly comes from Lemma \ref{psterr} replacing $\bv$ by $\bue$ and $\bv_h$ by $\Seh\bue$. For the second estimate, we differentiate (\ref{psve}) with respect to the temporal variable $t$ and do similar set of analysis as Lemma \ref{psterr}. This completes the estimates.
\end{proof}
\noindent 
Recall that, we split $\bxi_h$ as follows: $ \bxi_h = \bzeta+\btheta $. 
Armed with the estimates of $\bzeta$ and $\bzeta_t,$ we now pursue the estimates of $\btheta$ in order to find the optimal error estimates of $\bxi_h$ in $L^{\infty}(L^2)$ and $L^{\infty}(H^1)$-norms. From (\ref{pbxi}) and (\ref{psve}), the equation in $\btheta$ turns out to be
\begin{equation}\label{pbth}
 (\btheta_t,\bphi_h)+\nu a_{\ve}(\btheta,\bphi_h) = -(\bzeta_t,\bphi_h),~~~~ \forall \bphi_h\in\bH_h.
\end{equation}
%--------------------------------------------------L \infty of p\bxi_h
\begin{lemma}\label{lipbxi}
Under the assumptions of Lemma \ref{paph}, there is a positive constant $C$ such that  $\bxi_h$ satisfies the following estimate for some finite time $t>0$,
\begin{equation*}%\label{lipbxi01}
 \|\bxi_h(t)\| + h\|\nabla \bxi_h(t)\| ~\le~ Ch^{m+1}t^{-\frac{m}{2}}.
\end{equation*}
\end{lemma}
\begin{proof}
Selecting  $\bphi_h=\sigma^{m}(t)\btheta$ with $\sigma^{m}(t)=e^{2\alpha t}\tau^m(t)$ in (\ref{pbth}), it now follows that
\begin{align}\label{lipbxi003}
 \frac{1}{2}\frac{d}{dt}\big(\sigma^{m}(t)\|\btheta\|^2\big)+\nu\sigma^{m}(t)\|A_{\ve h}^{\frac{1}{2}}  \btheta\|^2 =& -\sigma^{m}(t)(\bzeta_t,\btheta)+\frac{1}{2}\sigma^{m}_{t}(t)\|\btheta\|^2.
\end{align}
A use of the Cauchy-Schwarz inequality with an appropriate use of the Young's inequality with $\sigma^{m}_{t}(t)\le (2\alpha+m)\sigma^{m-1}(t)$ in (\ref{lipbxi003}) yields
 \begin{equation}\label{lipbxi004}
 \frac{1}{2}\frac{d}{dt}\big(\sigma^{m}(t)\|\btheta\|^2\big)+\nu\sigma^{m}(t)\|A_{\ve h}^{\frac{1}{2}} \btheta\|^2  \le \frac{1}{2}\sigma^{m+1}(t)\|\bzeta_t\|^2 +\frac{1}{2}\big(2\alpha+m+1\big) \sigma^{m-1}(t)\|\btheta\|^2.
 \end{equation}
On integrating (\ref{lipbxi004}) with respect to time from $0$ to $t$  and write $\btheta=\bxi_h -\bzeta$ to find  that
\begin{align*} %\label{lipbxi005}
\sigma^{m}(t)\|\btheta(t)\|^2 +\nu \int_0^t \sigma^{m}(s)\|A_{\ve h}^{\frac{1}{2}}\btheta(s)\|^2\,ds 
\le C\int_0^t \sigma^{m+1}(s)\|\bzeta_s(s)\|^2\,ds  + C\int_0^t  \sigma^{m-1}(s) \|\btheta(s)\|^2 ds \nonumber \\
\le C\int_0^t \sigma^{m+1}(s)\|\bzeta_s(s)\|^2\,ds  + C\int_0^t  \sigma^{m-1}(s)(\|\bxi_h(s)\|^2+\|\bzeta(s)\|^2)~ ds.
\end{align*}
Use the Lemmas \ref{l2epbxi} and \ref{psvest0} to conclude
\begin{align*} %\label{theta}
\sigma^{m}(t)\|\btheta(t)\|^2+\int_0^t \sigma^{m}(s)\|A_{\ve h}^{\frac{1}{2}}\btheta(s)\|^2\,ds  \le  Ch^{2m+2}\Big(\int_0^t  \sigma^{m+1}(s)\|A_{\ve}^{\frac{m+1}{2}}\bues(s)\|^2 ds \nonumber\\
 + \int_0^t  \sigma^{m-1}(s)\|A_{\ve}^{\frac{m+1}{2}}\bue(s)\|^2 ds\Big).
\end{align*}
Finally, a use of Lemma \ref{psvest0} with Lemma \ref{Aep}, triangle inequality and the inverse hypothesis completes the rest of the proof.
\end{proof}

\noindent
Recall that
$$ \eve := \bue - \bueh = (\bue - \btveh) + ( \btveh - \bueh)= \bxi_h + \bta_h. $$
Since all the required estimates of $\bxi_h$ are obtained, it is now enough to estimate $\bta_h$. 
%
%------------------------------------------------------L^2 Penalized Error Estimate
%
\begin{lemma}\label{l2perr}
 Suppose the assumptions of Lemma \ref{paph} hold. Let $\bueh(t)$ be a solution of (\ref{dwfpomu}) with initial condition
$\bueh(0)= P_h\bu_{\ve 0}$.  Then, there exists a positive constant $C$ such that the following holds:
\begin{equation*} %\label{l2perr01}
 \int_0^t \sigma^{m-1}(s)\|\eve(s)\|^2~ds \le K(t)h^{2m+2}\int_0^t \sigma^{m-1}(s) \|A_{\ve}^{\frac{m+1}{2}}\hbue(s)\|^2 ds,\quad m\ge 1,
\end{equation*}
where, $K(t)=Ce^{Ct}$.
\end{lemma}
\begin{proof}
 In view of the Lemma \ref{l2epbxi}, we need to prove only the estimate of $\bta_h.$ From (\ref{dwfpomu}) and (\ref{lpo}), the equation in $\bta_h$ becomes
\begin{equation}\label{pbta}
 (\bta_{ht},\bphi_h)+\nu a_{\ve}(\bta_h,\bphi_h) =\tb(\bueh,\bueh,\bphi_h)-\tb(\bue,\bue,\bphi_h),~~~\bphi_h\in\bH_h.
\end{equation}
Choose $\bphi_h=A_{\ve h}^{-1}\bta_h$ to obtain
\begin{align}\label{l2perr001}
 \frac{1}{2}\frac{d}{dt}\|A_{\ve h}^{-\frac{1}{2}}\bta_h\|^2 +\nu\|\bta_h\|^2 = \Lambda_h(A_{\ve h}^{-1}\bta_h),
\end{align}
where
\begin{align*} %\label{nonlinear}
\Lambda_h(\bphi_h)= \tb(\bueh,\bueh,\bphi_h)-\tb(\bue,\bue,\bphi_h)
 = -\tb(\eve,\bueh,\bphi_h)-\tb(\bue,\eve,\bphi_h).
\end{align*}
Use Lemma \ref{trilinear} with Lemma \ref{Aeph} and write $\eve=\bxi_h+\bta_h$, we easily estimate $\Lambda_h$ as similar to \cite[(4.9)-(4.14)]{GD15} as
\begin{align}\label{l2perr003}
 |\Lambda_h(A_{\ve h}^{-1}\bta_h)| \le \rho\|\bta_h\|^2
 +C(\rho)\Big(\|\nabla\bueh\|^2+\|\bueh\|\|\nabla\bueh\|+\|\bue\|\|\nabla \bue\|\Big)\|\bxi_h\|^2 \nonumber\\
 +C(\rho)\|A_{\ve h}^{-\frac{1}{2}}\bta_h\|^2\Big(\|\nabla\bueh\|^4+\|\bueh\|^2\|\nabla
 \bueh\|^2 +\|\bue\|^2\|\nabla\bue\|^2\Big).
\end{align}
With $\rho=\nu/2,$ we obtain using (\ref{l2perr001}) and Lemma \ref{paph}
\begin{equation}\label{l2perr002}
 \frac{d}{dt}\|A_{\ve h}^{-\frac{1}{2}}\bta_h\|^2+\nu\|\bta_h\|^2 \le C(\nu)\Big(\|\bxi_h\|^2 +\|A_{\ve h}^{-\frac{1}{2}}\bta_h\|^2\Big).
\end{equation}
Now, multiply by $e^{2\alpha t}$ and apply the Gronwall's Lemma \cite{S95} to conclude
\begin{equation*}
 e^{2\alpha t}\|A_{\ve h}^{-\frac{1}{2}}\bta_h(t)\|^2+\nu\int_0^t e^{2\alpha s} \|\bta_h(s)\|^2 ds \le C e^{Ct}\int_0^t e^{2\alpha s}\|\bxi_h(s)\|^2 ds.
\end{equation*}
An application of the triangular inequality with Lemma \ref{l2epbxi} concludes the proof for $m=1$. For $m>1$, multiply (\ref{l2perr002}) by $\sigma^{m-1}(t)$ and using the fact $\frac{d}{dt}\sigma^{m-1}(t)\le C\sigma^{m-2}(t)$ and (\ref{l2perr003}) with $\rho=\nu/2$, we obtain
\begin{equation*}
 \frac{d}{dt}(\sigma^{m-1}(t)\|A_{\ve h}^{-\frac{1}{2}}\bta_h\|^2) +\nu\sigma^{m-1}(t)\|\bta_h\|^2 
  \le C\sigma^{m-2}(t)\|A_{\ve h}^{-\frac{1}{2}}\bta_h\|^2 + C \sigma^{m-1}(t)(\|\bxi_h\|^2 +\|A_{\ve h}^{-\frac{1}{2}}\bta_h\|^2).
\end{equation*}
Finally, arguing in exactly the same way as in the Lemma \ref{l2epbxi} and applying the Gronwall's Lemma \cite{S95}, the triangular inequality with Lemma \ref{l2epbxi}, we complete the  rest of the proof.
\end{proof}

\begin{lemma}\label{essbta}
Under the  assumption of Lemma \ref{l2perr}, the following holds:
\begin{equation*}
 \|\bta_h(t)\|+h\|\nabla\bta_h(t)\| \le K(t) h^{m+1}t^{-\frac{m}{2}}.
\end{equation*}
\end{lemma}

\begin{proof}
Choose $\bphi_h=\sigma^{m}(t)\bta_h$ in (\ref{pbta}) to obtain
\begin{equation*}
 \frac{1}{2}\frac{d}{dt}(\sigma^{m}(t)\|\bta_h\|^2)+\nu\sigma^{m}(t)\|A_{\ve h}^{\frac{1}{2}}\bta_h\|^2 
  = \frac{1}{2}\sigma^{m}_{t}(t)\|\bta_h\|^2 +\sigma^{m}(t)\Lambda_h(\bta_h).
\end{equation*}
We proceed same as  \cite[(4.16)-(4.20)]{GD15} and integrate to find that
\begin{align*}
 \sigma^{m}(t)\|\bta_h\|^2+& \nu\int_0^t \sigma^{m}(s)\|A_{\ve h}^{\frac{1}{2}}\bta_h(s)\|^2\le \frac{1}{2}(2\alpha+m) \int_0^t \sigma^{m-1}(s) \|\bta_h(s)\|^2ds  \\
    +& C\int_0^t \tau(s)\big(\|\nabla\bue(s)\|\|\Delta\bue(s)\|+\|\nabla\bueh(s)\|
 \|\Delta_h\bueh(s)\|\big) \sigma^{m-1}(s)\|\bta_h(s)\|^2.
\end{align*}
Apply  Lemmas \ref{pap}, \ref{paph}  and \ref{l2perr}, and then multiply the resulting inequality by $e^{-2\alpha t}$ to arrive at
$$ \tau^m(t)\|\bta_h\|^2+\nu e^{-2\alpha t} \int_0^t \sigma^{m}(s)\|A_{\ve h}^{\frac{1}{2}}\bta_h\|^2~ds ~\le  Ch^{2m+2}. $$
Since $\bta_h\in\bH_h$, we use inverse hypothesis to obtain an estimate for $\|\nabla\bta_h\|$, which completes the rest of the proof.
\end{proof}

\subsection{Error Estimates for the Pressure}

Subtract the second equation of (\ref{wfpom}) from the second equation of (\ref{dwfpom}) and obtain
\begin{equation}\label{peer}
(\pe-\peh,\chi_h)=   \frac{\nu}{\ve}(\nabla\cdot\eve, \chi_h).
\end{equation}
Choose $\chi_h=\peh-j_h\pe=e_p-(\pe-j_h\pe)$ with $e_p=\pe-\peh$ in (\ref{peer}) to find that
\begin{align}
\|e_p\|^2 &= (e_p,\pe-j_h\pe) +  \frac{\nu}{\ve}(\nabla\cdot\eve, e_p)- \frac{\nu}{\ve}(\nabla\cdot\eve, \pe-j_h\pe) \nonumber\\
&\le Ch^{2m}\|\pe\|_{m}^2+ \frac{C}{\ve^2}\|\nabla\cdot\eve\|^2 + \frac{1}{2}\|e_p\|^2.
\end{align}
If we use the bound $\|\nabla\cdot\bphi\| \le C\|\nabla\bphi\|$, then we observe that the error bound depends on $1/\ve$. We therefore find an alternate bound for divergent of error instead of replacing by gradient of velocity. Replacing $\bv$ by $\bue$ and $\bv_h$ by $\Seh\bue$ in (\ref{sterr11}), we can say $\|\nabla\cdot(\bue-\Seh\bue)\|\le C\sqrt{\ve}h^mt^{-\frac{m}{2}}.$ Since $\eve=(\bue-\Seh\bue)+(\Seh\bue-\bueh)$, hence, $ \frac{1}{\ve}\|\nabla\cdot\eve\|$ will depend on $1/\sqrt\ve$, and so will $e_p$. 
It is clear that the error bound for the pressure always depends on $1/\ve$ or $1/\sqrt\ve$ if we find it directly using velocity error. But, if we choose the finite element spaces $\bH_h$ and  $L_h$ in such a way that satisfy the discrete inf-sup condition ${\bf (B2^\prime)}$, then we can find the $\ve$-uniform pressure error estimate as given below. 

\noindent
First, we split $e_p$ as
\begin{equation*} %\label{p1.1}
\|e_p\|=\|\pe-\peh\| \le \|\pe-j_h\pe\| + \|j_h\pe-\peh\|. 
\end{equation*}
From ${\bf (B2^\prime)}$, we observe that
\begin{align}
\|j_h\pe-\peh\| &\le C \sup_{\bphi_h\in\bH_h\setminus\{0\}}\bigg\{\frac{|(j_h\pe-\peh, \nabla\cdot\bphi_h)|}{\|\nabla\bphi_h\|}\bigg\} \nonumber \\
    &\le  C\bigg(\|j_h\pe-\pe\|+ \sup_{\bphi_h\in\bH_h\setminus \{0\}}\bigg\{\frac{|(\pe-\peh, \nabla\cdot\bphi_h)|}{\|\nabla\bphi_h\|}\bigg\}\bigg). \label{p2}
\end{align}
The first term on right hand side of $(\ref{p2})$ can be estimated by using ${\bf (B1)}$ and for the second term, we subtract  (\ref{dwfpom}) from (\ref{wfpom}) to obtain
\begin{align}
  (e_p,\nabla\cdot\bphi_h)=(\evet,\bphi_h)+\nu a(\eve,\bphi_h)  -\Lambda_h(\bphi) ~~~\forall\bphi_h\in \bH_h. \label{p3}
\end{align}
A use of Lemma \ref{trilinear} shows
\begin{align}
|\Lambda_h(\bphi_h)| \le C(\|\nabla\bue\|+\|\nabla\bueh\|)\|\nabla\eve\|\|\nabla\bphi_h\|.  \label{p4}
\end{align} 
Now, apply the Cauchy-Schwarz inequality (\ref{p3}) and use (\ref{p4}) to arrive at
\begin{align}
  (e_p,\nabla\cdot\bphi_h)\le &\Big[C\|\evet\|_{-1;h}+ C \|\nabla\eve\| \Big]\|\nabla\bphi_h\|   \label{p5}
\end{align}
where, 
\begin{equation}\label{negnorm}
\|\evet\|_{-1;h} = \sup\Big\{\frac{<\evet,\bphi_h>}{\|\nabla\bphi_h\|}:\bphi_h\in\bH_h,\bphi_h\neq0\Big\}.
\end{equation}
Since all the estimate on right hand side in (\ref{p5}) are known except $\|\evet\|_{-1;h}$. Also $\|\evet\|_{-1;h} \le \|\evet\|_{-1}$, so now we derive $\|\evet\|_{-1}$. 

\begin{lemma}\label{pensdplm}%%penalty semi disc pressure lemma
The error $\eve$ satisfies for $0<t<T$
\[\|\evet\|_{-1}\le K(t) \Big(h^{m}(\|\buet\|_{m-1}+ \|\nabla\eve\|\Big).\]
\end{lemma}
\begin{proof}
For any $\bpsi\in \bH_0^1$, use the orthogonal projection $P_h:\bL^2\to \bH_h$, we obtain using (\ref{ppeesd}) with $\bphi_h=P_h\bpsi$
\begin{align}\label{pensdplm1}
(\evet,\bpsi)&=(\evet,\bpsi-P_h\bpsi) + (\evet,P_h\bpsi) \nonumber \\
    &= (\evet,\bpsi-P_h\bpsi) - \nu a_\ve(\eve,P_h\bpsi) - \Lambda_h(P_h\bpsi). 
\end{align}
Using approximation property of $P_h$, we find that
\begin{equation}
(\evet,\bpsi-P_h\bpsi)=(\buet-P_h\buet,\bpsi-P_h\bpsi) \le Ch^{m}\|\buet\|_{m-1}\|\nabla\bpsi\|. \label{pensdplm2}
\end{equation}
Also, using Lemma \ref{trilinear} with boundedness of $\bue$ and $\bueh$ to bound
\begin{equation}
\Lambda_h(P_h\bpsi) \le C(\|\nabla\bue\|+\|\nabla\bueh\|)\|\nabla\eve\|\|\nabla\bpsi\| \le C\|\nabla\eve\|\|\nabla\bpsi\|. \label{pensdplm3}
\end{equation}
Now substitute (\ref{pensdplm2})-(\ref{pensdplm3}) in (\ref{pensdplm1})  to obtain
\begin{align*}
(\evet,\bpsi) &\le C\Big(h^{m}(\|\buet\|_{m-1}+ \|\nabla\eve\|\Big)\|\nabla\bpsi\|.  %\label{sdplm4}
\end{align*}
and therefore,
\begin{align*} %\label{pensdplm5}
\|\evet\|_{-1} &\le \sup\Big\{\frac{<\evet,\bphi>}{\|\nabla\bphi\|}:\bphi\in\bH_0^1,\bphi\neq0\Big\} \nonumber \\
      & \le C\Big(h^{m}(\|\buet\|_{m-1}+ \|\nabla\eve\|\Big). 
\end{align*}
This completes the rest of the proof.
\end{proof}
The use of (\ref{p2})-(\ref{p5}) and Lemma \ref{pensdplm} and with the help of Lemma \ref{pap} and Theorem \ref{perrest} yields the following.
\begin{lemma}\label{thmp}
Under the hypothesis of the Lemma $\ref{l2epbxi}$, there exists a positive constant $C$ define as Lemma $\ref{l2epbxi}$, such that, for all $t>0$, it holds:
\[  \|(\pe-\peh)(t)\|  \le K(t) h^{m}t^{-\frac{m}{2}}.\]
\end{lemma}
\noindent
\textit{Proof of the Theorem \ref{perrest}.} 
Combining the Lemmas  \ref{lipbxi}, \ref{essbta} and  \ref{thmp}, we complete the proof the the Theorem \ref{perrest}.
\begin{remark}
It is noted that, there is no restriction on choosing finite element spaces for the velocity error bounds. But for the $\ve$-uniform pressure error bound we have to choose a proper finite element spaces which satisfy the discrete inf-sup condition. In \cite{HLY06}, it is seen that the improper choice of finite element spaces like ($P1-P0$) makes the error bounds dependent of $1/\sqrt\ve$ in the context of the steady state Navier-Stokes equations. 
\end{remark}

%%%% Uniform estimates
%
\subsection{Uniform Estimates}

We note here that the error estimates obtained in Theorem 4.1 are exponentially dependent on time.
In order to find the uniform estimates of $\eve$, we need to obtain the uniform estimates of $\bta_h$ only using (\ref{unique}), since, the estimates of $\bxi_h$ are uniform in time (see, Lemma \ref{lipbxi}). 
\begin{theorem}\label{uniest}
Under the assumptions of Theorem \ref{perrest} and the uniqueness condition (\ref{unique}), there exists a positive constant $C$ such that 
\begin{equation*}%\label{perrest0}
 \|(\bue-\bueh)(t)\|+h(\|\nabla(\bue-\bueh)(t)\|+\|(\pe-\peh)(t)\|) \le C h^{m+1}t^{-\frac{m}{2}},
\end{equation*}
\end{theorem}
\begin{proof}
Choose $\bphi_h=\bta_h$ in (\ref{pbta}) to obtain 
\begin{equation}\label{unbta1}
 \frac{1}{2}\frac{d}{dt}(\|\bta_h\|^2) + \nu \|A_{\ve h}^{\frac{1}{2}}\bta_h\|^2 = \Lambda_h(\bta_h).
\end{equation}
We use Lemma \ref{trilinear}, \ref{Aeph} and the uniqueness condition (\ref{unique}) with Lemma \ref{pap}, \ref{paph}, \ref{lipbxi} to bound $\Lambda_h(\bta_h)$ as
\begin{align*}
|\Lambda_h(\bta_h)| &= \tb(\bxi_h,\bueh,\bta_h)+\tb(\bue,\bxi_h,\bta_h)+\tb(\bta_h,\bueh,\bta_h)\nonumber\\
& \le C(\|\nabla\bue\|^{\frac{1}{2}}\|\Delta\bue\|^{\frac{1}{2}}+\|\nabla\bueh\|^{\frac{1}{2}}\|\Delta\bueh\|^{\frac{1}{2}})\|\bxi_h\|\|\nabla\bta_h\| + N \|\nabla\bueh\|\|\nabla\bta_h\|^2 \nonumber\\
& \le Ch^{m+1}\tau^{-\frac{2m+1}{4}}\|\nabla\bta_h\|  + N  \|\nabla\bueh\|\|\nabla\bta_h\|^2.
\end{align*}
Use the above estimate in (\ref{unbta1}) and multiply by $e^{2\alpha t}$ then integrate with respect to time  to arrive at
\begin{align*}
 e^{2\alpha t} \|\bta_h(t)\|^2  + 2\int_0^t \big(\nu-N \|\nabla\bueh\|\big)e^{2\alpha s} \|\nabla\bta_h(s)\|^2 ds  + \frac{2\nu}{\ve}\int_0^t  e^{2\alpha s}\|\nabla\cdot\bta_h(s)\|^2 ds \\
  \le \|\bta_h(0)\|^2  + 2\alpha \int_0^t e^{2\alpha s}\|\bta_h(s)\|^2 ds + Ch^{m+1} \int_0^t \tau^{-\frac{2m+1}{4}}(s)e^{2\alpha s}\|\nabla\bta_h(s)\| ds.
\end{align*}
Now, multiply by $e^{-2\alpha t}$ and take $t\to \infty$ and use 
$\overline{\lim}_{t\to\infty}\|\nabla\bueh\| \le \frac{1}{\nu}\|\f\|_{L^\infty(0,\infty;\bH^{-1}(\Omega))}$  and $\overline{\lim}_{t\to\infty} \tau^{-\frac{2m+1}{4}}(t) = 1$ to find that
\begin{equation*}
\big(\nu-\frac{N}{\nu}\|\f\|_{L^\infty(0,\infty;\bH^{-1}(\Omega))}\big)\overline{\lim_{t\to\infty}} \|\nabla \bta_h(t)\|^2  \le  Ch^{m+1}  \overline{\lim_{t\to\infty}}\|\nabla\bta_h(t)\|.
\end{equation*}
From (\ref{unique}), we have $ N\nu^{-2} \|\f\|_{L^\infty(0,\infty;\bH^{-1}(\Omega))}<1$. Then, we conclude that
\begin{equation*}
\overline{\lim_{t\to\infty}} \|\bta_h(t)\|  \le \overline{\lim_{t\to\infty}} \|\nabla \bta_h(t)\|  \le  Ch^{m+1}.
\end{equation*}
With this and  Lemma \ref{lipbxi}, we finally obtain
\begin{equation*}
\overline{\lim_{t\to\infty}} \|\eve(t)\|  \le  Ch^{m+1}t^{-\frac{m}{2}}.
\end{equation*}
Using the above estimate, we easily derive the estimates of $\|\nabla\eve(t)\|$ and $\|(\pe-\peh)(t)\|$, which completes the rest of the proof.
\end{proof}

\section{Backward Euler Method}
\se

This section focuses on a completely discrete method based on the backward Euler method. 
For time discretization, let $k,~0<k<1,$ be the time step and $t_n=nk,~n\ge 0.$ We define for a sequence $\{\bphi^n\}_{n
\ge 0}\subset\bH_h,$ the backward difference quotient
$$ \pt\bphi^n=\frac{1}{k}(\bphi^n-\bphi^{n-1}). $$
For any continuous function $\bv(t)$ we set $\bv_n=\bv(t_n)$. 
We describe below the backward Euler scheme for the penalized semidiscrete NSEs  (\ref{dwfpom}): Find $\{\bUe^n\}_{1\le n\le N}\in\bH_h$ and $\{\Pe^n\}_{1\le n\ge N}\in L_h$ as solutions of the recursive nonlinear algebraic equations ($1\le n\le N$) :
%-----------------------------------------------fully dicsrete (backward Euler) in H
\begin{equation}\left\{\begin{array}{l}\label{fdbeh}
 (\pt\bUe^n,\bphi_h)+\nu a(\bUe^n,\bphi_h) = (\Pe^n,\nabla\cdot \bphi_h)
 +(\f^n,\bphi_h)-\tb(\bUe^n,\bUe^n,\bphi_h),\forall \bphi_h\in\bH_h\\
 \nu(\nabla\cdot\bUe^n,\chi_h)+\ve(\Pe^n,\chi_h) =0~~~\forall \ \chi_h \in L_h,~~~n\ge 0,
\end{array}\right.
\end{equation}
where, $\bUe^0=P_h\bu_{\ve 0}$.
As earlier we can get rid of the fully discrete penalized pressure term
%---------------------------------------------fully dicsrete (backward Euler) in J
\begin{equation}\label{fdbej}
 (\pt\bUe^n,\bphi_h)+\nu a_{\ve}(\bUe^n,\bphi_h)= (\f^n,\bphi_h)
 -\tb(\bUe^n,\bUe^n,\bphi_h)~~~\forall\bphi_h\in\bH_h.
\end{equation}
Using variant of Brouwer fixed point theorem and standard uniqueness arguments, it is easy to show that the discrete problem (\ref{fdbej}) is well-posed. For a proof, we refer to the Appendix. 
Below, we prove {\it a priori} bounds for the discrete solutions $\{\bU^n\}_{1\le n\le N}.$
We note here that since the bounds proved below are independent of $n,~1\le n\le N$, these bounds are uniform in time, that is, estimates are still valid as the final time $t_N\to +\infty$.

\begin{lemma}\label{stb}
Let $\alpha_0>0$ be such that for $0<\alpha<\min\{\alpha_0, \frac{\nu\lambda_1}{2c_0^2}\}$,
$$1+(\frac{\nu\lambda_1}{2c_0^2}) k \ge e^{2\alpha k}.$$
Then, the discrete solution $\bUe^n,~1\le n\le N$ of (\ref{fdbej}) satisfies the following estimates:
\begin{align}
\|\bUe^n\|^2+  \frac{\nu}{2} e^{-2\alpha t_{n}}k\sum_{i=1}^n e^{2\alpha t_i}\|A_{\ve h}^{\frac{1}{2}}\bUe^i\|^2 &\le C, \label{stb1a} \\
\tau_n^{(r-1)}\|A_{\ve h}^{\frac{r}{2}}\bUe^n\|^2+  \nu e^{-2\alpha t_n}k\sum_{i=1}^n e^{2\alpha t_i}\|A_{\ve h}\bUe^i\|^2 &\le C,\quad r\in\{1,2\}, \label{stb1b}
\end{align}
where $C$ depends on the given data. 
\end{lemma}

\begin{proof}
For $n=i$, we  substitute $\bphi_h=\bUe^i$ in (\ref{fdbej}) and use $(\pt\bUe^i,\bUe^i) \geq \frac{1}{2}\pt\|\bUe^i\|^2$ and the Poincar\'e inequality and the Cauchy-Schwarz inequality with Lemma \ref{Aeph} $(\|\bUe^i\|\le \frac{1}{\sqrt\lambda_1}\|\nabla\bUe^i\|\le \frac{c_0}{\sqrt\lambda_1}\|A_{\ve h}^{\frac{1}{2}}\bUe^i\|)$ to obtain
\begin{align}
\frac{1}{2}\pt\|\bUe^i\|^2 +  \frac{3\nu}{4} \|A_{\ve h}^{\frac{1}{2}}\bUe^i\|^2  \le  \frac{c_0^2}{\nu\lambda_1} \|\f^i\|^2. \label{lemmabe4.6}
\end{align}
Note that, the nonlinear term $\tb(\bUe^i,\bUe^i,\bUe^i)=0$. Hence multiply (\ref{lemmabe4.6}) by $ke^{2\alpha t_i}$ and sum over $i=1$ to $n$ and use the fact 
\begin{align*}
k\sum_{i=1}^n e^{2\alpha t_i}\pt\|\bUe^i\|^2 
&= e^{2\alpha t_n}\|\bUe^n\|^2 -  \|\bUe^0\|^2 - k\sum_{i=1}^{n-1} \Big(\frac{e^{2\alpha k}-1}{k}\Big) e^{2\alpha t_i}\|\bUe^i\|^2\\
& \ge e^{2\alpha t_n}\|\bUe^n\|^2 -  \|\bUe^0\|^2 - k\sum_{i=1}^{n} c_0^2\Big(\frac{e^{2\alpha k}-1}{k\lambda_1}\Big) e^{2\alpha t_i}\|A_{\ve h}^{\frac{1}{2}}\bUe^i\|^2,
\end{align*}
to conclude that
\begin{align*}
e^{2\alpha t_n}\|\bUe^n\|^2 + \Big(\frac{3\nu}{2}- c_0^2\big(\frac{ e^{2\alpha k}-1}{k\lambda_1}\big)\Big) k \sum_{i=1}^{n}e^{2\alpha t_i}\|A_{\ve h}^{\frac{1}{2}}\bUe^i\|^2  \le \|\bUe^0\|^2 + \frac{2c_0^2\|\f\|_{\infty}^2}{\nu\lambda_1} k \sum_{i=1}^{n}e^{2\alpha  t_i},
\end{align*}
where $\|\f\|_{\infty}=\|\f\|_{L^\infty(\mathbb{R}_+;\bL^2)}$. 
With $0<\alpha<\min\{\alpha_0, \frac{\nu\lambda_1}{2c_0^2}\}$, which guarantees that $\nu \ge c_0^2\big(\frac{ e^{2\alpha k}-1}{k\lambda_1}\big)$. 
Now, we multiply both sides by $e^{-2\alpha t_n}$ to conclude (\ref{stb1a}).\\
\noindent For (\ref{stb1b}), we substitute $\bphi_h=A_{\ve h}\bUe^i$ in (\ref{fdbej}) to obtain
\begin{align}\label{lemma4.61}
\frac{1}{2}\pt\|A_{\ve h}^{\frac{1}{2}}\bUe^i\|^2 + \nu \|A_{\ve h}\bUe^i\|^2  \le  (\f^i,A_{\ve h}\bUe^i)- \tb(\bUe^i,\bUe^i,A_{\ve h}\bUe^i).
\end{align}
Now a use of the Cauchy-Schwarz inequality and the Young's inequality yields
\begin{align}\label{lemmabe4.62}
|(\f^i,A_{\ve h}\bUe^i)- \tb(\bUe^i,\bUe^i,A_{\ve h}\bUe^i)| &\le \|\f^i\|\|A_{\ve h}\bUe^i\|+ C \|\bUe^i\|^{\frac{1}{2}}\|A_{\ve h}^{\frac{1}{2}}\bUe^i\| \|A_{\ve h}\bUe^i\|^{3/2} \nonumber\\
&\le C(\|\f^i\|^2+  \|\bUe^i\|^{2}\|A_{\ve h}^{\frac{1}{2}}\bUe^i\|^4) + \frac{\nu}{2} \|A_{\ve h}\bUe^i\|^{2}.
\end{align}
Insert (\ref{lemmabe4.62}) in (\ref{lemma4.61}) to obtain
\begin{align}\label{lemma4.63}
\pt\|A_{\ve h}^{\frac{1}{2}}\bUe^i\|^2 + \nu \|A_{\ve h}\bUe^i\|^2  \le C\big(  \|\f^i\|^2  
+    \|\bUe^i\|^{2}\|A_{\ve h}^{\frac{1}{2}}\bUe^i\|^4 \big).
\end{align}
Now, choose $\bphi_h=\pt\bUe^n$ in (\ref{fdbej}) with $n=i$ and use the fact $(\pt\bphi_h^n,\bphi_h^n)=\frac{1}{2}\pt\|\bphi_h^n\|^2+\frac{k}{2}\|\pt\bphi_h^n\|^2$ to obtain
\begin{align}\label{lemma_dt1}
\|\pt\bUe^i\|^2 + \frac{\nu}{2}\pt\|A_{\ve h}^{\frac{1}{2}}\bUe^i\|^2 + \frac{k\nu}{2}\|\pt A_{\ve h}^{\frac{1}{2}}\bUe^i\|^2 \le (\f^i, \pt\bUe^i) - \tb(\bUe^i,\bUe^i,\pt\bUe^i). 
\end{align}
Using (\ref{tbp1}), the last term on the right hand side of (\ref{lemma_dt1}) can be written as
\begin{align*}
|\tb(\bUe^i,\bUe^i,\pt\bUe^i)| 
=\frac{1}{k}|\tb(\bUe^i,\bUe^i,\bUe^i-\bUe^{i-1})| = \frac{1}{k}|\tb(\bUe^i,\bUe^{i},\bUe^{i-1})|=\frac{1}{k}|\tb(\bUe^i,\bUe^{i-1},\bUe^{i})|.
\end{align*}
A use of Lemma \ref{trilinear} with \ref{Aeph}, the Cauchy-Schwarz inequality and the Young's inequality yields
\begin{align}\label{lemma_dtin}
\frac{1}{k}|\tb(\bUe^i,\bUe^{i-1},\bUe^{i})|
& \le \frac{C}{k} \|\bUe^i\|^{\frac{1}{2}}\|\nabla\bUe^i\|^{\frac{1}{2}} \|\nabla\bUe^{i-1}\| \|\bUe^{i}\|^{\frac{1}{2}} \|\nabla\bUe^i\|^{\frac{1}{2}}  \nonumber\\
& \le \frac{C}{k} \|\bUe^i\| \|A_{\ve h}^{\frac{1}{2}}\bUe^i\|  \|A_{\ve h}^{\frac{1}{2}}\bUe^{i-1}\|  \nonumber\\
& \le \frac{\nu}{4k}\|A_{\ve h}^{\frac{1}{2}}\bUe^i\|^2 + \frac{C}{k\nu} \|\bUe^i\|^2 \|A_{\ve h}^{\frac{1}{2}}\bUe^{i-1}\|^{2}.
\end{align}
We use the inequality (\ref{lemma_dtin}) with $\|\bUe^i\|\le C$ in (\ref{lemma_dt1}) and use the fact $\pt\|A_{\ve h}^{\frac{1}{2}}\bUe^i\|^2= \frac{1}{k}(\|A_{\ve h}^{\frac{1}{2}}\bUe^i\|^2-\|A_{\ve h}^{\frac{1}{2}}\bUe^{i-1}\|^2)$, then we obtain
\begin{align*} %\label{lemma_dtin1}
\|\pt\bUe^i\|^2 + \frac{\nu}{4k}\|A_{\ve h}^{\frac{1}{2}}\bUe^i\|^2 + \frac{k\nu}{4}\|\pt A_{\ve h}^{\frac{1}{2}}\bUe^i\|^2 \le \frac{c_0^2}{k\nu}\|\f^i\|^2+\frac{C}{k}\| A_{\ve h}^{\frac{1}{2}}\bUe^{i-1}\|^2. 
\end{align*}
Now, drop the first and the last terms on the left hand side and multiply the resulting equation by $4k/\nu$ to arrive at
\begin{align}\label{lemma_dt2}
 \|A_{\ve h}^{\frac{1}{2}}\bUe^i\|^2 \le C \|A_{\ve h}^{\frac{1}{2}}\bUe^{i-1}\|^2 + C \|\f^i\|^2. 
\end{align}
A use of (\ref{lemma_dt2}) in the last term on the right hand side of (\ref{lemma4.63}) leads to
\begin{align}\label{lemma_dt3}
\pt\|A_{\ve h}^{\frac{1}{2}}\bUe^i\|^2 + \nu  \|A_{\ve h}\bUe^i\|^2  \le g^{i-1} \|A_{\ve h}^{\frac{1}{2}}\bUe^{i-1}\|^2 + h^{i-1},
\end{align}
where, $g^{i-1}= C \|\bUe^i\|^{2}\|A_{\ve h}^{\frac{1}{2}}\bUe^i\|^2$ and $h^{i-1}=C\big(1+\|\bUe^i\|^{2}\|A_{\ve h}^{\frac{1}{2}}\bUe^i\|^2\big) \|\f^i\|^2$. We now use the discrete uniform Gronwall's lemma (Lemma \ref{dugl}) to derive 
\begin{align}\label{lemma_dt4}
\|A_{\ve h}^{\frac{1}{2}}\bUe^n\|^2 \le C, \quad \forall n \ge N_1+1,
\end{align}
where $N_1\le N$. For $1\le n\le N_1$, we use the classical discrete Gronwall's lemma \cite{HR90,S90} to obtain
\begin{align}\label{lemma_dt41}
\|A_{\ve h}^{\frac{1}{2}}\bUe^n\|^2 \le C, \quad \text{for}~ 1\le n \le N_1.
\end{align}
We multiply (\ref{lemma_dt3}) by $e^{2\alpha t_i}$ and sum from $i=1$ to $n$  and use (\ref{stb1a}), (\ref{lemma_dt4}) and (\ref{lemma_dt41}).  Finally, we multiply the resulting equation by $e^{-2\alpha t_n}$ to conclude (\ref{stb1b}), when $r=1$.\\
For $r=2$, we need a couple of intermediate estimates. For this, first we use Lemma \ref{trilinear} with Lemma \ref{Aeph} and the Cauchy-Schwarz inequality to bound the terms on the right hand side of (\ref{lemma_dt1}) as 
\begin{align}\label{lemma_dt42}
|(\f^i, \pt\bUe^i) - \tb(\bUe^i,\bUe^i,\pt\bUe^i)| \le  C(\|\f^i\|^2 + \|A_{\ve h}^{\frac{1}{2}}\bUe^i\|^{2}\|A_{\ve h}\bUe^i\|^2)+ \frac{1}{2}\|\pt\bUe^i\|^2.
\end{align} 
After using (\ref{lemma_dt42}) in (\ref{lemma_dt1}), we multiply the resulting equation by $ke^{2\alpha t_i}$ and take sum from $i=1$ to $n$. Then a use of (\ref{stb1b}) with $r=1$ gives
\begin{align}\label{lemma_dt43}
\nu \|A_{\ve h}^{\frac{1}{2}}\bUe^n\|^2 + e^{-2\alpha t_n}k\sum_{i=1}^n e^{2\alpha t_i}\|\pt\bUe^i\|^2  \le C.
\end{align}
For the second intermediate estimate, we consider (\ref{fdbej}) for $n=i$ and for $n=i-1$. Then  subtract them and divide by $k$ to obtain
\begin{equation}\label{dfdbej}
 (\pt(\pt\bUe^i),\bphi_h)+\nu a_{\ve}(\pt\bUe^i,\bphi_h)= (\pt\f^i,\bphi_h)
 -\tb(\pt\bUe^i,\bUe^i,\bphi_h)-\tb(\bUe^{i-1},\pt\bUe^i,\bphi_h)~~~\forall\bphi_h\in\bH_h.
\end{equation}
Now, choose $\bphi_h=\pt\bUe^i$ in (\ref{dfdbej}) and use Lemma \ref{trilinear} with Lemma \ref{Aeph} and the Cauchy-Schwarz inequality to find
\begin{align}\label{lemma_dt44}
 \frac{1}{2}\pt\|\pt\bUe^i\|^2 + \nu \|A_{\ve h}^{\frac{1}{2}}\pt\bUe^i\|^2 
 &\le  \|\pt\f^i\|\|\pt\bUe^i\|
 + C \|\pt\bUe^i\|^{\frac{1}{2}}\|\nabla\pt\bUe^i\|^{\frac{1}{2}}\|\nabla\bUe^i\|\|\pt\bUe^i\|^{\frac{1}{2}}\|\nabla\pt\bUe^i\|^{\frac{1}{2}} \nonumber\\
 & \le C \Big( \|\pt\f^i\|^2
 + \|\pt\bUe^i\|^2 \|A_{\ve h}^{\frac{1}{2}}\bUe^i\|^2\Big) + \frac{\nu}{2}\|A_{\ve h}^{\frac{1}{2}}\pt\bUe^i\|^2.
\end{align}
Multiply (\ref{lemma_dt44}) by $k\tau_i e^{2\alpha t_i}$ and sum from $i=1$ to $n$ and use (\ref{lemma_dt43}) to deduce 
\begin{equation}\label{lemma_dt5}
\tau_n\|\pt\bUe^n\|^2 \le C.
\end{equation} 
Now, choose $\bphi_h=A_{\ve h}\bUe^n$ in (\ref{fdbej}) and use (\ref{lemmabe4.62}) to rewrite it as
\begin{align*}
 \nu \|A_{\ve h}\bUe^n\|^2  
 & =  (\f^n,A_{\ve h}\bUe^n)- \tb(\bUe^n,\bUe^n,A_{\ve h}\bUe^n)- (\pt \bUe^n, A_{\ve h}\bUe^n)\nonumber\\
 & \le C\big(\|\f^n\|+  \|\bUe^n\|\|A_{\ve h}^{\frac{1}{2}}\bUe^n\|^2 + \|\pt\bUe^n\|\big) \|A_{\ve h}\bUe^n\|.
\end{align*}
Multiply by $\tau_n^{\frac{1}{2}}$ and use (\ref{stb1a}), (\ref{lemma_dt4}) and (\ref{lemma_dt5}) to conclude (\ref{stb1b}), when $r=2$, which completes the proof.
\end{proof}
\begin{remark}  %\label{alpha_0}
In Lemma \ref{stb}, such a choice of $\alpha_0>0$ is possible by choosing $\alpha_0 < \frac{\log(1+\frac{\nu\lambda_1}{2c_0^2}k)}{k}$. Note that for large $k>0$, $\alpha_0$ is small but as $k\to 0$, $\frac{\log(1+\frac{\nu\lambda_1}{2c_0^2}k)}{k} \to \frac{\nu\lambda_1}{2c_0^2}$. Thus, with $0<\alpha < \min\{\alpha_0,\frac{\nu\lambda_1}{2c_0^2}\}$, the result in Lemma \ref{stb} is valid.
\end{remark}

\subsection{Fully Discrete Error Estimate}

Define $\bueh(t_n)=\bueh^n$ and $\eve^n=\bUe^n-\bueh^n$. Consider (\ref{dwfpomu}) at $t=t_n$ and subtract from (\ref{fdbej}) to arrive at
\begin{align}\label{bee1}
(\pt\eve^n,\bphi_h)+\nu a_{\ve}(\eve^n,\bphi_h) = R_h^n(\bphi_h)+\Lambda_h^n(\bphi_h),
\end{align}
where
\begin{align}
R_h^n(\bphi_h)&=(\bueht^{n},\bphi_h)-(\pt\bueh^n,\bphi_h)=(\bueht^n,\bphi_h)-\frac{1}{k}\int_{t_{n-1}}^{t_n}(\buehs,\bphi_h)~ds \label{RH}\nonumber\\
&=\frac{1}{k}\int_{t_{n-1}}^{t_n}(t-t_{n-1})(\buehss,\bphi_h)~ds, \\
\Lambda_h^n(\bphi_h)&= -\tb(\bueh^n,\bueh^n,\bphi_h)-\tb(\bUe^n,\bUe^n,\bphi_h)= -\tb(\eve^n,\bueh^n,\bphi_h)-\tb(\bUe^n,\eve^n,\bphi_h),\label{LH}
\end{align}

The optimal error estimates for velocity and pressure for time discretization only, have already been proved in \cite{S95} when the initial data $\bu_{\ve 0}\in \bH^2\cap\bH_0^1$ and under the assumption of sufficiently small $k$ so that $Ck<1$. In our case, $\bu_{\ve 0}\in \bH_0^1$ and without smallness assumption on $k$, we have proved the following optimal error estimates.

\begin{lemma}  \label{eeL1}
Assume that $({\bf A1})$,$({\bf A2})$, $({\bf B1})$ and $({\bf B2})$ hold true. Then, under the assumption of Lemma \ref{stb} there exists some positive constant $C$, that depends on the given data, the following holds
\begin{align*}
\|A_{\ve h}^{j/2}\eve^n\|^2 + k e^{-2\alpha t_n}\sum_{i=1}^n e^{2\alpha t_i} \|A_{\ve h}^{(j+1)/2}\eve^i\|^2 
\leq K_n k^{1-j},~~ j= -1,0,1,
\end{align*}
where $K_n=Ce^{Ct_n}$.
\end{lemma}
\begin{proof}
For $j=0$, take $n=i$ and $\bphi_h= \eve^i$ in (\ref{bee1}) to arrive at
\begin{align}
\pt \|\eve^i\|^2 + 2\nu \|A_{\ve h}^{\frac{1}{2}}\eve^i\|^2  \le  2R_h^i( \eve^i) +2\Lambda_h^i( \eve^i). \label{eeL1.11}
\end{align}
Then multiply (\ref{eeL1.11}) by $k e^{2\alpha t_i}$ and take sum from $i=1$ to $n$ and use the fact
\begin{align}\label{eeL1_dt}
k\sum_{i=1}^n e^{2\alpha t_i}\pt\|\eve^i\|^2 
&\geq e^{2\alpha t_n}\|\eve^n\|^2 -  \sum_{i=1}^{n-1} e^{2\alpha t_{i}}(e^{2\alpha k}-1)\|\eve^{i}\|^2 \nonumber\\
&\geq e^{2\alpha t_n}\|\eve^n\|^2 - c_0^2\Big(\frac{e^{2\alpha k}-1}{k\lambda_1}\Big) k\sum_{i=1}^{n-1} e^{2\alpha t_{i}}\|A_{\ve h}^{\frac{1}{2}}\eve^{i}\|^2 ,
\end{align}
to find that 
\begin{align}
e^{2\alpha t_n}\|\eve^n\|^2 &+ \Big(2\nu-c_0^2\big(\frac{e^{2\alpha k}-1}{k\lambda_1}\big)\Big) k\sum_{i=1}^n e^{2\alpha t_i}\|A_{\ve h}^{\frac{1}{2}}\eve^i\|^2 
\le  2k\sum_{i=1}^{n} e^{2\alpha t_i}\big(R_h^i(\eve^i)+\Lambda_h^i(\eve^i)\big). \label{eeL1.21}
\end{align}
A use of the Cauchy-Schwarz's inequality with (\ref{RH}) and $t-t_{i-1}\le t_i, t\in [t_{i-1},t_i]$ yields
\begin{align}\label{eeL1.31}
k\sum_{i=1}^{n}e^{2\alpha t_i} R_h^i(\eve^i) &\le k\sum_{i=1}^n  e^{2\alpha t_i} \Big(\frac{1}{k}\int_{t_{i-1}}^{t_i} (s-t_{i-1})\|A_{\ve h}^{-\frac{1}{2}}\buehss\|~ ds\Big) \|A_{\ve h}^{\frac{1}{2}}\eve^i\| \nonumber\\
	&\le   \frac{1}{\nu k}\sum_{i=1}^n   \bigg(\int_{t_{i-1}}^{t_i}(s-t_{i-1})~ ds\bigg)\bigg( \int_{t_{i-1}}^{t_i} e^{2\alpha t_i} (s-t_{i-1}) \|A_{\ve h}^{-\frac{1}{2}}\buehss\|^2 ds\bigg) \nonumber\\ &~~~~~+ \frac{\nu}{4} k\sum_{i=1}^n e^{2\alpha t_i} \|A_{\ve h}^{\frac{1}{2}}\eve^i\|^2 \nonumber\\
	&\le C k \int_0^{t_n} \sigma(s) \|A_{\ve h}^{-\frac{1}{2}}\buehss\|^2 ds + \frac{\nu}{4} k\sum_{i=1}^n e^{2\alpha t_i} \|A_{\ve h}^{\frac{1}{2}}\eve^i\|^2.
\end{align}
From Lemma \ref{Aeph} and \ref{trilinear}, we find that
\begin{align}
k\sum_{i=1}^{n} e^{2\alpha t_i} \Lambda_h^i(\eve^i) &\le C k\sum_{i=1}^{n} e^{2\alpha t_i}  \|\nabla\bu_h^i\|^2  \|\eve^i\|^2 + \frac{\nu}{4} k\sum_{i=1}^{n} e^{2\alpha t_i} \|A_{\ve h}^{\frac{1}{2}}\eve^i\|^2 \nonumber\\
	&\le C k\sum_{i=1}^{n} e^{2\alpha t_i} \|\eve^i\|^2 + \frac{\nu}{4} k\sum_{i=1}^{n} e^{2\alpha t_i} \|A_{\ve h}^{\frac{1}{2}}\eve^i\|^2\nonumber\\
	&\le C k e^{2\alpha t_n} \|\eve^n\|^2 + C k\sum_{i=1}^{n-1} e^{2\alpha t_i} \|\eve^i\|^2 + \frac{\nu}{4} k\sum_{i=1}^{n} e^{2\alpha t_i} \|A_{\ve h}^{\frac{1}{2}}\eve^i\|^2.\label{eeL1.41}
\end{align}
Inserting (\ref{eeL1.31})-(\ref{eeL1.41}) in (\ref{eeL1.21}) with the boundedness of $\|\eve^n\| \le \|\bue^n\|+ \|\bUe^n\|\le C$ and Lemma \ref{paph}, we conclude that
\begin{align*}
e^{2\alpha t_n}\|\eve^n\|^2 &+ \Big(\nu-c_0^2\big(\frac{e^{2\alpha k}-1}{k\lambda_1}\big)\Big) k\sum_{i=1}^n e^{2\alpha t_i}\|A_{\ve h}^{\frac{1}{2}}\eve^i\|^2 
\le C k e^{2\alpha t_n} + C k\sum_{i=1}^{n-1} e^{2\alpha t_i} \|\eve^i\|^2.
\end{align*}
With $0<\alpha<\min\{\alpha_0, \frac{\nu\lambda_1}{2c_0^2}\}$, we have $\big(\nu-c_0^2\big(\frac{e^{2\alpha k}-1}{k\lambda_1}\big)\big)>0$. 
Finally, a use of the discrete Gronwall's lemma \cite{HR90,S90} completes the rest of proof for the case $j=0$.
For the case $j=-1$, take $n=i$ and $\bphi_h= A_{\ve h}^{-1}\eve^i$ in (\ref{bee1}) and multiply the resulting equation by $k e^{2\alpha t_i}$ and take sum from $i=1$ to $n$ and use the similar fact (\ref{eeL1_dt}) to obtain
\begin{align}
e^{2\alpha t_n}\|A_{\ve h}^{-\frac{1}{2}}\eve^n\|^2 + \Big(2\nu-c_0^2\big(\frac{e^{2\alpha k}-1}{k\lambda_1}\big)\Big) k\sum_{i=1}^n e^{2\alpha t_i}\|\eve^i\|^2 
\le  2k\sum_{i=1}^{n} e^{2\alpha t_i}\big(R_h^i(A_{\ve h}^{-1}\eve^i)+\Lambda_h^i(A_{\ve h}^{-1}\eve^i)\big). \label{eeL1.2}
\end{align}
A use of the Cauchy-Schwarz's inequality with (\ref{RH}) yields
\begin{align}\label{eeL1.3}
k\sum_{i=1}^{n}e^{2\alpha t_i} R_h^i(A_{\ve h}^{-1}\eve^i) &\le k\sum_{i=1}^n  e^{2\alpha t_i} \Big(\frac{1}{k}\int_{t_{i-1}}^{t_i} (s-t_{i-1})\|A_{\ve h}^{-1}\buehss\|~ ds\Big) \|\eve^i\| \nonumber\\
	&\le C  k\sum_{i=1}^n   \bigg(\int_{t_{i-1}}^{t_i}  ds\bigg)\bigg( \int_{t_{i-1}}^{t_i} e^{2\alpha t_i} \|A_{\ve h}^{-1}\buehss\|^2 ds\bigg) + \frac{\nu}{4} k\sum_{i=1}^n e^{2\alpha t_i} \|\eve^i\|^2 \nonumber\\
	&\le C k^2 \int_0^{t_n} e^{2\alpha s} \|A_{\ve h}^{-1}\buehss\|^2 ds + \frac{\nu}{4} k\sum_{i=1}^n e^{2\alpha t_i} \|\eve^i\|^2.
\end{align}
From Lemma  \ref{Aeph} and \ref{trilinear}, it follows that
\begin{align}
k\sum_{i=1}^{n} e^{2\alpha t_i} \Lambda_h^i(A_{\ve h}^{-1}\eve^i) &\le C k\sum_{i=1}^{n} e^{2\alpha t_i} (\|\nabla\bu_h^i\|^2+\|\nabla\bU_h^i\|^2)\|A_{\ve h}^{-\frac{1}{2}}\eve^i\|^2 + \frac{\nu}{4} k\sum_{i=1}^{n} e^{2\alpha t_i} \|\eve^i\|^2 \nonumber\\
	&\le C k\sum_{i=1}^{n} e^{2\alpha t_i} \|A_{\ve h}^{-\frac{1}{2}}\eve^i\|^2 + \frac{\nu}{4} k\sum_{i=1}^{n} e^{2\alpha t_i} \|\eve^i\|^2 \nonumber\\
	&\le C k  e^{2\alpha t_n} \|A_{\ve h}^{-\frac{1}{2}}\eve^n\|^2 + C k\sum_{i=1}^{n-1} e^{2\alpha t_i} \|A_{\ve h}^{-\frac{1}{2}}\eve^i\|^2 + \frac{\nu}{4} k\sum_{i=1}^{n} e^{2\alpha t_i} \|\eve^i\|^2.\label{eeL1.4}
\end{align}
Using (\ref{eeL1.3})-(\ref{eeL1.4}) in (\ref{eeL1.2}) with $\|A_{\ve h}^{-\frac{1}{2}}\eve^n\|^2 \le \|\eve^n\|_{-1}^2 \le C \|\eve^n\|^2 \le Ck$, we obtain
\begin{align*}
e^{2\alpha t_n}\|A_{\ve h}^{-\frac{1}{2}}\eve^n\|^2 +  \Big(\nu-c_0^2\big(\frac{e^{2\alpha k}-1}{k\lambda_1}\big)\Big) k\sum_{i=1}^n e^{2\alpha t_i}\|\eve^i\|^2 
\le Ck^2 e^{2\alpha t_n} + C \sum_{i=1}^{n-1} e^{2\alpha t_i}\|A_{\ve h}^{-\frac{1}{2}}\eve^i\|^2.
\end{align*}
Then, a use of the discrete Gronwall's lemma \cite{HR90,S90} concludes the result for the case $j=-1$. 
For the last one, that is, $j=1$, we choose $\bphi_h=A_{\ve h} \eve^i $ and follow a similar analysis as above to complete the rest of the proof.
\end{proof}

\begin{lemma}  \label{eeL2}
Let the assumption of Lemma \ref{eeL1} be satisfied. Then, for some positive constant $C$, that depends on $T$, there holds
\begin{align*}
\tau_n\|\eve^n\|^2 + k e^{-2\alpha t_n}\sum_{i=1}^n \sigma_i \|A_{\ve h}^{\frac{1}{2}}\eve^i\|^2 \leq K_n k^2,
\end{align*}
where $\sigma_{i}=\tau_i e^{2\alpha t_i}$ and $\tau_i=\tau(t_i)=\min\{1,t_i\}$.
\end{lemma}

\begin{proof}
Take $n=i$ and $\bphi_h=\sigma_i \eve^i$ in (\ref{bee1}) to arrive at
\begin{align}
 \sigma_i\pt\|\eve^i\|^2 +2\nu\sigma_i\|A_{\ve h}^{\frac{1}{2}}\eve^i\|^2  \le  2R_h^i(\sigma_i \eve^i) +2\Lambda_h^i(\sigma_i \eve^i). \label{beb1}
\end{align}
Now multiply (\ref{beb1}) by $k$ and take sum from $i=1$ to $n$ and use the fact $$k\sum_{i=1}^n\sigma_i\pt\|\eve^i\|^2 \geq \sigma_n\|\eve^n\|^2 - k\sum_{i=1}^{n-1} e^{2\alpha t_{i}}\|\eve^{i}\|^2,$$ to obtain
\begin{align}
 \sigma_n\|\eve^n\|^2 &+ 2\nu k\sum_{i=1}^n \sigma_i\|A_{\ve h}^{\frac{1}{2}}\eve^i\|^2 \le k\sum_{i=1}^{n-1} e^{2\alpha t_{i}}\|\eve^{i}\|^2 + 2k\sum_{i=1}^{n} \sigma_i\big(R_h^i(\eve^i)+\Lambda_h^i(\eve^i)\big). \label{beb2}
\end{align}
A use of the Cauchy-Schwarz inequality with the Young's inequality, (\ref{RH}) and Lemma \ref{Aep1} yields
\begin{align}\label{beb3}
k\sum_{i=1}^{n} \sigma_i R_h^i(\eve^i) &\le k\sum_{i=1}^n  \sigma_i \Big(\frac{1}{k}\int_{t_{i-1}}^{t_i} (s-t_{i-1})\|A_{\ve h}^{-\frac{1}{2}}\buehss\|~ ds\Big)~ \|A_{\ve h}^{\frac{1}{2}}\eve^i\| \nonumber\\
	&\le C  k\sum_{i=1}^n   \bigg(\int_{t_{i-1}}^{t_i}  ds\bigg)\bigg( \int_{t_{i-1}}^{t_i} \sigma_i(s)\|A_{\ve h}^{-\frac{1}{2}}\buehss\|^2 ds\bigg) + \frac{\nu}{4} k\sum_{i=1}^n \sigma_i \|A_{\ve h}^{\frac{1}{2}}\eve^i\|^2 \nonumber\\
	&\le C k^2 \int_0^{t_n} \sigma(s)\|A_{\ve h}^{-\frac{1}{2}}\buehss\|^2 ds + \frac{\nu}{4} k\sum_{i=1}^n \sigma_i \|A_{\ve h}^{\frac{1}{2}}\eve^i\|^2.
\end{align}
From the Lemmas \ref{trilinear}, \ref{Aeph} and \ref{eeL1}, we deduce that
\begin{align}
k\sum_{i=1}^{n} \sigma_i \Lambda_h(\eve^i) &\le k\sum_{i=1}^{n} \sigma_i \big(\|\Delta_h\bueh^i\|^2+\|\Delta_h\bUe^i\|^2\big)\|\eve^i\|^2 + \frac{\nu}{4} k\sum_{i=1}^{n} \sigma_i \|\nabla\eve^i\|^2 \nonumber\\
	&\le C k\sum_{i=1}^{n} e^{2\alpha t_i}\|\eve^i\|^2 + \frac{\nu}{4} k\sum_{i=1}^{n} \sigma_i \|A_{\ve h}^{\frac{1}{2}}\eve^i\|^2.\label{beb4}
\end{align}
Using (\ref{beb3})-(\ref{beb4}) in (\ref{beb2}) with Lemma \ref{eeL1}, we conclude the rest of the result.
\end{proof}
\begin{lemma}  \label{eeL3}
Let the assumption of Lemma \ref{eeL1} be satisfied. Then, for some positive constant $C$, that depends on $T$, there holds
\begin{align*}
\tau_n^2\|A_{\ve h}^{\frac{1}{2}}\eve^n\|^2 + k e^{-2\alpha t_n}\sum_{i=1}^n \sigma_{i}^{2} \|A_{\ve h}\eve^i\|^2 \leq K_n k^2,
\end{align*}
where $\sigma_{i}^{2}=\tau_i^2 e^{2\alpha t_i}$.
\end{lemma}
\begin{proof}
The proof is quite similar to the proof of the previous lemma. So we only give a sketch of the proof. 
Choose $\bphi_h=\sigma_{i}^{2}A_{\ve h} \eve^n$ with $n=i$ in (\ref{bee1}) to arrive at
\begin{align*}
 \sigma_i^2\pt\|A_{\ve h}^{\frac{1}{2}}\eve^i\|^2 +2\nu\sigma_i^2\|A_{\ve h}\eve^i\|^2  \le  2R_h^i(\sigma_i^2 A_{\ve h}\eve^i) +2\Lambda_h^i(\sigma_i^2 A_{\ve h}\eve^i). %\label{bebh1}
\end{align*}
Multiplying by $k$ and summing from $i=1$ to $n$, we obtain
\begin{align}
 \sigma_n^2\|A_{\ve h}^{\frac{1}{2}}\eve^n\|^2 &+ 2\nu k\sum_{i=1}^n \sigma_i^2\|A_{\ve h}\eve^i\|^2 \le k\sum_{i=1}^{n-1} \sigma_i\|A_{\ve h}^{\frac{1}{2}}\eve^{i}\|^2 + 2k\sum_{i=1}^{n} \sigma_i^2\big(R_h^i(A_{\ve h}\eve^i)+\Lambda_h^i(A_{\ve h}\eve^i)\big). \label{bebh2}
\end{align}
We apply the Cauchy-Schwarz inequality and the Young's inequality with (\ref{RH}) and Lemma \ref{Aep1} to bound
\begin{align}\label{bebh3}
k\sum_{i=1}^{n} \sigma_i^2 R_h^i(A_{\ve h}\eve^i) &\le k\sum_{i=1}^n  \sigma_i^2 \Big(\frac{1}{k}\int_{t_{i-1}}^{t_i} (s-t_{i-1})\|\buehss\|~ ds\Big)~ \|A_{\ve h}\eve^i\| \nonumber\\
	&\le C  k\sum_{i=1}^n   \bigg(\int_{t_{i-1}}^{t_i}  ds\bigg)\bigg( \int_{t_{i-1}}^{t_i} \sigma_i^2(s)\|\buehss\|^2 ds\bigg) + \frac{\nu}{4} k\sum_{i=1}^n \sigma_i^2 \|A_{\ve h}\eve^i\|^2 \nonumber\\
	&\le C k^2 \int_0^{t_n} \sigma^2(s)\|\buehss\|^2 ds + \frac{\nu}{4} k\sum_{i=1}^n \sigma_i^2 \|A_{\ve h}\eve^i\|^2.
\end{align}
With the help of the Lemmas \ref{trilinear}, \ref{Aeph} and \ref{eeL1}, we can bound the nonlinear terms as
\begin{align}
k\sum_{i=1}^{n} \sigma_i^2 \Lambda_h(A_{\ve h}\eve^i) &\le k\sum_{i=1}^{n} \tau_i \big(\|\Delta_h\bueh^i\|^2+\|\Delta_h\bUe^i\|^2\big) \sigma_i\|A_{\ve h}^{\frac{1}{2}}\eve^i\|^2 + \frac{\nu}{4} k\sum_{i=1}^{n} \sigma_i^2 \|A_{\ve h}\eve^i\|^2 \nonumber\\
	&\le C k\sum_{i=1}^{n} \sigma_i\|A_{\ve h}^{\frac{1}{2}}\eve^i\|^2 + \frac{\nu}{4} k\sum_{i=1}^{n} \sigma_i^2 \|A_{\ve h}\eve^i\|^2.\label{bebh4}
\end{align}
Using (\ref{bebh3})-(\ref{bebh4}) in (\ref{bebh2}) with the Lemmas \ref{eeL1} and \ref{eeL3}, we conclude the rest of the proof.
\end{proof}

\noindent
We now also derive the error estimate for the pressure term. In fact,  similar to the semidiscrete pressure error estimate, we can easily prove that $\tau_n\|\pt\eve^n\|\le K_n k$. Now using this and the available estimates for $\eve^n$, we can easily prove the following lemma:
\begin{lemma}  \label{eeL4}
Let the assumption of Lemma \ref{eeL1} be satisfied. Then, for some positive constant $C$, that depends on $T$, there holds
\begin{align*}
 \tau_n \|\Pe^n - \peh(t_n)\| \leq K_n k.
\end{align*}
\end{lemma}
\begin{proof}
From (\ref{dwfpom}) and (\ref{fdbeh}), we can write the pressure error equation as
\begin{align*}
(P_{\ve}^n - \peh^n,\nabla\cdot\bphi_h) = (\pt\eve^n,\bphi_h)+\nu a(\eve^n,\bphi_h) - R_h^n(\bphi_h) - \Lambda_h^n(\bphi_h),
\end{align*}
where $R_h^n$ and $\Lambda_h^n$ are defined by (\ref{RH}) and (\ref{LH}), respectively.
A use of Lemma \ref{trilinear} gives
\begin{align*}
(P_{\ve}^n - \peh^n,\nabla\cdot\bphi_h) = \Big(\|\pt\eve^n\|_{-1,h} + \nu \|\nabla\eve^n\| + C(\|\nabla\bueh^n\|+\|\nabla\bUe^n\|)\|\nabla\eve^n\| \nonumber \\ + \frac{1}{k}\int_{t_{n-1}}^{t_n}(t-t_{n-1})\|\buehss\|_{-1}~ds\Big)  \|\nabla\bphi_h\|,
\end{align*}
where $\|\cdot\|_{-1,h}$ is defined in (\ref{negnorm}) and clearly $\|\cdot\|_{-1,h} \le \|\cdot\|_{-1 }\le C \|\cdot\|$.
Finally, we use the Lemmas \ref{paph}, \ref{stb} and \ref{eeL3} to complete the rest of the proof.
\end{proof}

\begin{lemma}  \label{fduniform}
Let the assumption of Lemma \ref{eeL1} be satisfied. Then, 
\begin{align}\label{fduni}
\sqrt\tau_n\|\eve^n\| + \tau_n\|\nabla\eve^n\| +  \tau_n\|\Pe^n - \peh(t_n)\| \leq C k,
\end{align}
where $C$ depends exponentially on time. Under the uniqueness condition (\ref{unique}), the above estimate (\ref{fduni}) holds uniformly in time.
\end{lemma}
\begin{proof}
A use of Lemma \ref{eeL2}, \ref{eeL3} and \ref{eeL4} with triangular inequality shows (\ref{fduni}). Since the constant $C$ depends exponentially in time due to the use of Lemma \ref{eeL1} which is not uniform in time.  For improving the estimates of Lemma \ref{eeL1} using (\ref{unique}),
we multiply (\ref{bee1}) by $k e^{2\alpha t_i}$ and take sum from $i=i_0+1$ to $n$ and use (\ref{eeL1_dt})
to obtain 
\begin{align}\label{fduniform1}
e^{2\alpha t_n}\|\eve^n\|^2 &+ \Big(2\nu-c_0^2\big(\frac{e^{2\alpha k}-1}{k\lambda_1}\big)\Big) k\sum_{i=i_0+1}^n e^{2\alpha t_i}\|A_{\ve h}^{\frac{1}{2}}\eve^i\|^2 \nonumber\\
&\le e^{2\alpha t_{i_0}}\|\eve^{i_0}\|^2 +  2k\sum_{i=i_0+1}^{n} e^{2\alpha t_i}\big(R_h^i(\eve^i)+\Lambda_h^i(\eve^i)\big).
\end{align}
From (\ref{eeL1.31}), it follows that
\begin{align}\label{fduniform2}
2k\sum_{i=i_0+1}^{n}e^{2\alpha t_i} R_h^i(\eve^i) \le C k e^{2\alpha t_n} + \frac{\nu}{2} k\sum_{i=i_0+1}^n e^{2\alpha t_i} \|A_{\ve h}^{\frac{1}{2}}\eve^i\|^2.
\end{align}   
Now, we bound the nonlinear terms using (\ref{unique}) as
\begin{align}
|\Lambda_h^i(\eve^i)| \le N \|\nabla\bu_h^i\| \|\nabla\eve^i\|^2.
\end{align}
From (\ref{pap001}), we can easily derive that $\overline{\lim}_{t\to\infty}\|\nabla\bue\| \le \nu^{-1}\|\f\|_{L^\infty(0,T;\bH^{-1})}$, which implies
\begin{equation}\label{fduniform3}
2k\sum_{i=i_0+1}^{n} e^{2\alpha t_i} \Lambda_h^i(\eve^i) \le 2N \nu^{-1}\|\f\|_{L^\infty(0,T;\bH^{-1})} ~ k\sum_{i=i_0+1}^{n} e^{2\alpha t_i} \|\nabla\eve^i\|^2.
\end{equation}
Inserting (\ref{fduniform2})-(\ref{fduniform3}) in (\ref{fduniform1}), we conclude that
\begin{align*}
e^{2\alpha t_n}\|\eve^n\|^2 &+ \Big(\nu-c_0^2\big(\frac{e^{2\alpha k}-1}{k\lambda_1}\big)\Big) k\sum_{i=i_0+1}^n e^{2\alpha t_i}\|A_{\ve h}^{\frac{1}{2}}\eve^i\|^2 \\
&\qquad + \Big( \frac{\nu}{2} - 2N \nu^{-1}\|\f\|_{L^\infty(0,T;\bH^{-1})}\Big) k\sum_{i=i_0+1}^n e^{2\alpha t_i}\|\nabla\eve^i\|^2 +  \frac{\nu}{2\ve} k\sum_{i=i_0+11}^n e^{2\alpha t_i}\|\nabla\cdot\eve^i\|^2 \\
&\le C k e^{2\alpha t_n} + e^{2\alpha t_{i_0}}\|\eve^{i_0}\|^2.
\end{align*}
With $0<\alpha<\min\{\alpha_0, \frac{\nu\lambda_1}{2c_0^2}\}$, we have $\big(\nu-c_0^2\big(\frac{e^{2\alpha k}-1}{k\lambda_1}\big)\big)>0$ and from (\ref{unique}), $\big( \frac{\nu}{2} - 2N \nu^{-1}\|\f\|_{L^\infty(0,T;\bH^{-1})}\big)\ge 0$. Then, we obtain
\begin{align}\label{fduniform4}
\|\eve^n\|^2 &+ e^{-2\alpha t_n} k\sum_{i=1}^n e^{2\alpha t_i}\|A_{\ve h}^{\frac{1}{2}}\eve^i\|^2 
\le C k.
\end{align}
Next, we choose $\bphi_h= A_{\ve h}^{-1}\eve^n$ in (\ref{bee1}) with $n=i$ and multiply the resulting equation by $k e^{2\alpha t_i}$ and take sum from $i=i_0+1$ to $n$. Then, arguing similar set of analysis as above one can obtain
\begin{align}
\|A_{\ve h}^{-\frac{1}{2}}\eve^n\|^2 + e^{-2\alpha t_n} k\sum_{i=1}^n e^{2\alpha t_i}\|\eve^i\|^2 
\le  Ck^2. \label{fduniform5}
\end{align}
Finally, a use of (\ref{fduniform4}) and (\ref{fduniform5}) instead of Lemma \ref{eeL1} in the proof of Lemma \ref{eeL2}, \ref{eeL3} and \ref{eeL4} complete the rest of proof.

\end{proof}
\noindent Finally, combining Theorem \ref{pperrest}, \ref{perrest}, and Lemma \ref{fduniform}, we conclude the following theorem:
\begin{theorem}\label{final}
Assume that $({\bf A1})$,$({\bf A2})$, $({\bf B1})$ and $({\bf B2})$ hold true. Then, for some positive constant $C$, that depends on $T$, there holds:
\begin{align*}
 \|\bu(t_n)-\bUe^n\| \le K_n\Big((\ve  + k) t^{-\frac{1}{2}} + h^{m+1}t^{-\frac{m}{2}}\Big),\\
 \|\nabla(\bu(t_n)-\bUe^n)\| \le K_n \Big( (\ve +k) t^{-1} + h^{m}t^{-\frac{m}{2}}\Big),\\
 \|p(t_n)-P_{\ve}^n\| \leq K_n \Big((\ve  + k) t^{-1} + h^{m}t^{-\frac{m}{2}}\Big),
\end{align*}
where the positive constant $K_n=Ce^{Ct_n}$ depends exponentially on time. 
The estimates are uniform in time under the uniqueness condition (\ref{unique}), that is, the constant $K_n$ becomes $C$.
\end{theorem}
\begin{remark}
Although we have discussed conforming finite element spaces on this article, but all the results remain valid for ($P_1^{NC}-P_0$) nonconforming elements. Therefore, the present analysis improves upon the results of Lu and Lin \cite{LL10} in the sense that optimal estimates in $L^2$-norm are obtained when initial data are in $\bH_0^1$. To this effect, a numerical experiment is presented in Section 6.
\end{remark}
\section{Numerical Experiments}
In this section, we present some numerical experiments that verify the results of previous section, mainly verify the order of convergence of the error estimates. 

We consider the NSEs in the domain $\Omega=[0,1]\times[0,1]$ subject to homogeneous Dirichlet boundary conditions. We approximate the equation using $P_2-P_1$, $P_3-P_2$ and $P_1^{NC}-P_0$ elements over a triangulation of $\Omega$. The domain is partitioned into triangles with size $h=2^{-i}, i=1,2,\dots,6$. To verify the theoretical result, we consider the following examples.
\begin{example} \label{ex1}
	For the experiment in 2D, we take the forcing term $f(x,y,t)$ such that the solution of the problem to be
	\begin{align*}
	u_1(x,y,t)&= 2e^{t} x^2(x-1)^2 y(y-1)(2y-1), \\
    u_2(x,y,t)&= -2e^{t} x(x-1)(2x-1) y^2(y-1)^2, \\
    p(x,y,t)&=   2e^t(x-y).
	\end{align*}
\end{example}
\begin{table}[h]  % P2P1
	\centering %
		\begin{tabular}{|c|c|c|c|c|c|c|c|}
		\hline
   h   & $\|u(t_n)-U^n\|_{L^2}$ & Rate & $\|u(t_n)-U^n\|_{H^1}$ & Rate & $\|p(t_n)-P^n\|_{L^2}$ & Rate  \\
   \hline
  1/2 	&  3.30633896e-03  &		 &  2.96918336e-02  & 		  &  3.29192462e-02 &		\\
  1/4   &  5.11077157e-04  & 2.6936  &  8.36078857e-03  & 1.8284  &  6.96272136e-03 &  2.2412  \\
  1/8   &  5.25170055e-05  & 3.2826  &  1.99271639e-03  & 2.0689  &  9.29102388e-04  & 2.9057  \\
  1/16  &  6.32080598e-06  & 3.0546  &  5.32350596e-04  & 1.9042  &  2.78763334e-04  & 1.7368  \\
  1/32  &  7.91350016e-07  & 2.9977  &  1.35504728e-04  & 1.9740  &  7.93302459e-05  & 1.8131   \\
  1/64  &  9.89942025e-08  & 2.9989  &  3.39036941e-05  & 1.9988  &  2.01622898e-05  & 1.9762   \\
   \hline   
	\end{tabular}
	\caption{Errors and convergence rates with $ \nu= 1, k = \ve= \mathcal{O}(h^3)$ at time $t=1$ using $P_2-P_1$ element}
	\label{p2p1}
\end{table}
\begin{table}[h]  %P3-P2
	\centering %
		\begin{tabular}{|c|c|c|c|c|c|c|c|}
		\hline
   h   & $\|u(t_n)-U^n\|_{L^2}$ & Rate & $\|u(t_n)-U^n\|_{H^1}$ & Rate & $\|p(t_n)-P^n\|_{L^2}$ & Rate  \\
   \hline 
  1/2   &  8.72519439e-04  &         &  9.64515599e-03  &         &  2.42266708e-02 &         \\
  1/4   &  7.86457181e-05  & 3.4717  &  2.53372020e-03  & 1.9285  &  3.77601493e-03  & 2.6816  \\
  1/8   &  5.40305188e-06  & 3.8635  &  3.35107263e-04  & 2.9156  &  4.04997117e-04  & 3.2209  \\
  1/16  &  3.65287201e-07  & 3.8866  &  4.34389795e-05  & 2.9476  &  3.73502090e-05  & 3.4387   \\
  1/32  &  2.40020778e-08  & 3.9278  &  5.55651552e-06  & 2.9667  &  3.34112960e-06  & 3.4827   \\
   \hline   
	\end{tabular}
	\caption{Errors and convergence rates with $ \nu= 1, k = \ve = \mathcal{O}(h^4)$ at time $t=1$ using $P_3-P_2$ element}
	\label{p3p2}
\end{table}
%	
%Results for P1NC-P0 element
\begin{table}[h]  % P1NC-P0 flage 1
	\centering %
		\begin{tabular}{|c|c|c|c|c|c|c|c|}
		\hline
   h   & $\|u(t_n)-U^n\|_{L^2}$ & Rate & $\|u(t_n)-U^n\|_{H^1}$ & Rate & $\|p(t_n)-P^n\|_{L^2}$ & Rate  \\
   \hline 
  1/2   &  1.55373715e-01  & 		 &  7.54243397e-01  &  		  &  9.17725216e-01  &      \\
  1/4   &  6.46328013e-02  & 1.2654  &  4.53947780e-01  & 0.7325  &  4.68679354e-01  & 0.9695       \\
  1/8   &  2.01782694e-02  & 1.6795  &  2.39739250e-01  & 0.9211  &  1.74839152e-01  & 1.4226  \\
  1/16  &  5.43929542e-03  & 1.8913  &  1.21753766e-01  & 0.9775  &  5.78140714e-02  & 1.5965  \\
  1/32  &  1.39082972e-03  & 1.9675  &  6.12053289e-02  & 0.9922  &  1.98791319e-02  & 1.5401   \\
  1/64  &  3.49954196e-04  & 1.9907  &  3.06603901e-02  & 0.9973  &  7.92197390e-03  & 1.3273   \\
   \hline   
	\end{tabular}
	\caption{Errors and convergence rates with $ \nu= 1, k = \ve= \mathcal{O}(h^2)$ at time $t=1$ using $P_1^{NC}-P_0$ element}
	\label{p1ncp0}
\end{table}
In Tables \ref{p2p1} and \ref{p3p2}, we present the numerical errors and convergence rates obtained on successive meshes for the penalty finite element method with the backward Euler scheme, applied to the system \eqref{om04}-\eqref{ibc04} using $P_m-P_{m-1}$ elements for $m=2,3$, respectively. The numerical analysis shows that the rate of convergence are $\mathcal{O}(h^{m+1})$ in $\bL^2$-norm and  $\mathcal{O}(h^m)$ in $\bH^1$-norm for the velocity and $\mathcal{O}(h^m)$ in $\bL^2$-norm for the pressure with the choice of $k=\ve=\mathcal{O}(h^{m+1})$ at the final time level, that is, when $t=1$ and take $\nu=1$. These results support the optimal convergence rates obtained in Theorem \ref{final}. 
In Table \ref{p1ncp0}, we present the numerical results for $P_1^{NC}-P_0$ element. It is observed in Table \ref{p1ncp0}
% as well as Fig \ref{figp1ncp0} 
that the rate of convergence in $\bL^2$-norm and $\bH^1$-norm of the velocity is 2 and 1, respectively. Moreover it is linear in pressure in $L^2$-norm.

\begin{example} \label{ex3d}
	For the experiment in 3D, we take the forcing term $f(x,y,z,t)$ such that the solution of the problem to be
	\begin{align*}
	u_1(x,y,z,t)&= \pi e^{t} \sin(\pi x)(\cos(\pi y)\sin(\pi z) - \cos(\pi z)\sin(\pi y)), \\
    u_2(x,y,z,t)&= \pi e^{t} \sin(\pi y)(\cos(\pi z) \sin(\pi x) - \cos(\pi x)\sin(\pi z)), \\
    u_3(x,y,z,t)&= \pi e^{t} \sin(\pi z)(\cos(\pi x)\sin(\pi y) - \cos(\pi y) \sin(\pi x)), \\
    p(x,y,z,t)&=   e^t (\sin(\pi x)\sin(\pi y)\sin(\pi z)-8.0/\pi^3)).
	\end{align*}
\end{example}
\begin{table}[h]  % P2P1 3d
	\centering %
		\begin{tabular}{|c|c|c|c|c|c|c|c|}
		\hline
   h   & $\|u(t_n)-U^n\|_{L^2}$ & Rate & $\|u(t_n)-U^n\|_{H^1}$ & Rate & $\|p(t_n)-P^n\|_{L^2}$ & Rate  \\
   \hline
  1/2 	&  4.39800806e-02  &		 &  4.97101276e-01  & 		  &  3.84531508e-01 &		\\
  1/4   &  4.80285161e-03  & 3.1949  &  9.96770050e-02  & 2.3182  &  6.41136630e-02 &  2.5844  \\
  1/8   &  4.08823054e-04  & 3.5543  &  1.92203894e-02  & 2.3746  &  1.56891870e-02  & 2.0308  \\
  1/16  &  5.02080312e-05  & 3.0254  &  4.63294928e-03  & 2.0526  &  3.90792118e-03  & 2.0053 \\
   \hline   
	\end{tabular}
	\caption{Errors and convergence rates with $ \nu= 1, k = \ve= \mathcal{O}(h^3)$ at time $t=0.1$ using $P_2-P_1$ element}
	\label{p2p13d}
\end{table}
\begin{table}[!h]  % P3P2 3d
	\centering %
		\begin{tabular}{|c|c|c|c|c|c|c|c|}
		\hline
   h   & $\|u(t_n)-U^n\|_{L^2}$ & Rate & $\|u(t_n)-U^n\|_{H^1}$ & Rate & $\|p(t_n)-P^n\|_{L^2}$ & Rate  \\
   \hline
  1/2 	&  1.45489050e-02  &		 &  3.13468857e-01  & 		  &  1.52818314e-01 &		\\
  1/4   &  8.52910943e-04  & 4.0923  &  3.87617913e-02  & 3.0156  &  4.07810264e-02 &  1.9058  \\
  1/8   &  5.90735824e-05  & 3.8518  &  5.43288281e-03  & 2.8348  &  3.73728889e-04  & 3.4478  \\
  1/16  &  3.78219433e-06  & 3.9652  &  7.03179542e-04  & 2.9497  &  4.21715883e-05  & 3.1476 \\
   \hline   
	\end{tabular}
	\caption{Errors and convergence rates with $ \nu= 1, k = \ve= \mathcal{O}(h^4)$ at time $t=0.1$ using $P_3-P_2$ element}
	\label{p3p23d}
\end{table}
We also present the numerical errors and convergence rates for Example \ref{ex3d} using $P_m-P_{m-1}$ elements for $m=2,3$, in Tables \ref{p2p13d} and \ref{p3p23d}, respectively. From the Tables \ref{p2p13d} and \ref{p3p23d}, we observe that the rate of convergence are $\mathcal{O}(h^{m+1})$ in $\bL^2$-norm and  $\mathcal{O}(h^m)$ in $\bH^1$-norm for the velocity and $\mathcal{O}(h^m)$ in $\bL^2$-norm for the pressure with the choice of $k=\ve=\mathcal{O}(h^{m+1})$ at the final time level, that is, when $t=0.1$ and take $\nu=1$. These results support the optimal convergence rates obtained in Theorem \ref{final}. 

\begin{figure}[h!]
\centering
\begin{tikzpicture}
\node at (0, -0.3) {$(0,0)$};
\node at (5, -0.3) {$(1,0)$};
\node at (0, 5.3) {$(0,1)$};
\node at (5, 5.3) {$(1,1)$};
\node at (2.5, -0.3) {$u_1=0,u_2=0$};
\node at (2.5, 5.3) {$u_1=1,u_2=0$};
\node at (-0.6, 2.6) {$u_1=0$};
\node at (-0.6, 2.3) {$u_2=0$};
\node at (5.6, 2.6) {$u_1=0$};
\node at (5.6, 2.3) {$u_2=0$};

 \draw (0,0) rectangle (5,5);
 % \draw[->,blue] (0.5,5.6) -- (4.5,5.6);

\node (A) at (0.5,5.6) {};
\node (B) at (4.5,5.6) {};
% arrows
\draw[->,line width=0.3mm] (A) -- (B);
\end{tikzpicture}
\caption{Domain $\Omega$ for lid-driven cavity flow.}
\label{lid}
\end{figure}
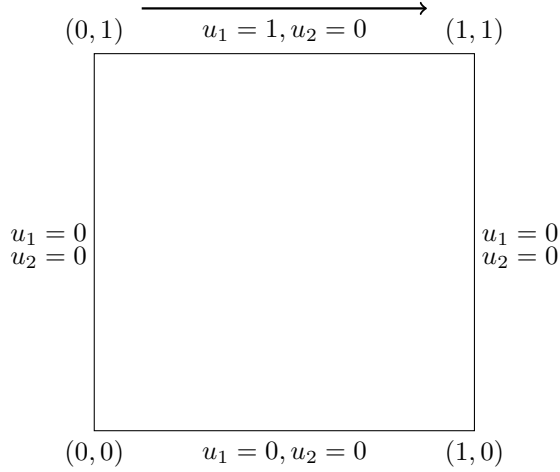

The next two examples are related to two-dimensional Benchmark problems.
\begin{example}\label{ex2}
In this example, we consider a benchmark problem related to a two-dimensional lid driven cavity flow on a unit square with zero body force. Also, no slip boundary condition are considered everywhere except the non  zero velocity $\bu=(1,0)^T$ on upper boundary, see Figure \ref{lid}.
\end{example}

\begin{figure}[h!] 
\centering
\includegraphics[scale=.55]{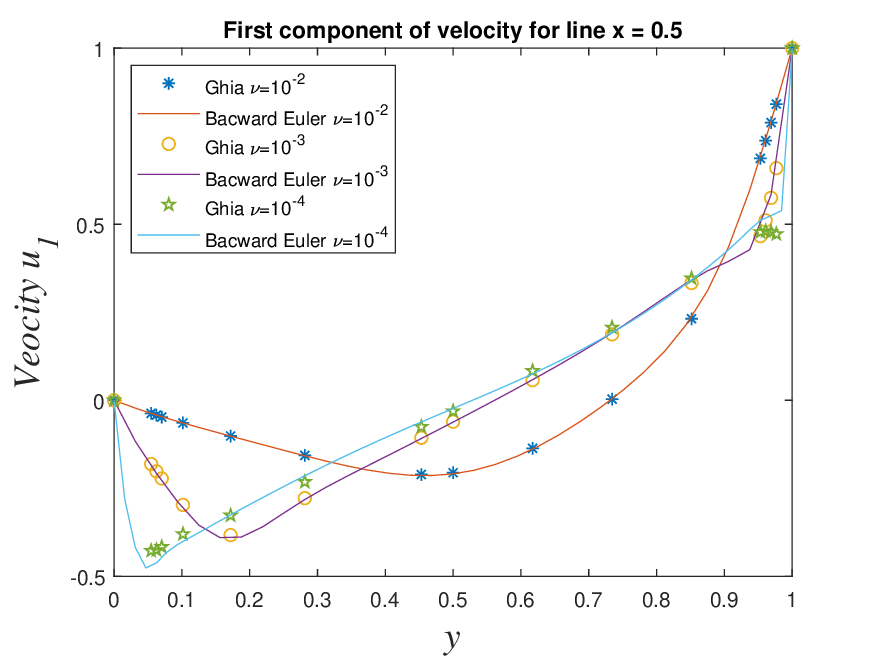}
\includegraphics[scale=.55]{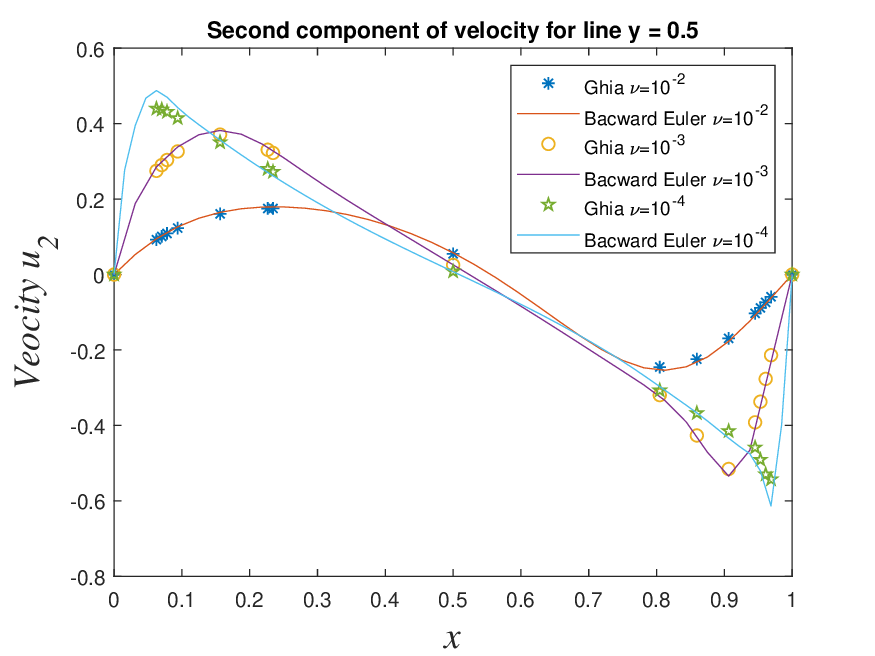}
\caption{Velocity components for Example \ref{ex2}.}
\label{fig7}
\end{figure}

For numerical experiments, we have chosen lines $(0.5,y)$ and $(x,0.5)$ and we plot the velocity profile with respect to these two lines. In Figure \ref{fig7}, we present the comparison between velocity obtained by penalty method and velocity obtained by Ghia et. al. \cite{GGS82} of NSEs for final time $t=75$, for $\nu=10^{-2}, 10^{-3}$ and $t=150$, for $\nu=10^{-4}$, respectively, with the choice of time step $k=0.01$. From the graphs, it is observed that the velocity profiles coincide with those of Ghia's results very well for a large time and that for $\nu$ small.
%
% Compare with Ghia

\begin{example}\label{ex2f}
Now, we also consider a well-known benchmark problem related to the two-dimensional flow around a cylinder with zero body forces \cite{J04}. The domain $\Omega$ is the channel of size $[0,2.2]\times[0,0.41]$ with a circle of diameter $0.1$ located at $(0.2,0.2)$ as shown in Figure \ref{fpc}. The whole boundary is divided in four parts; inflow boundary $\Gamma_{in} := \{x=0\}$, outflow boundary $\Gamma_{out} := \{x=2.2\}$, the remaining two wall $\Gamma_{wall} := \{y=0, y=0.41\}$ and the boundary of the circle $\Gamma_{cyl} := \{(x-0.2)^2+(y-0.2)^2=0.0025\}$. We consider the no-slip boundary at $\Gamma_{wall}$ and $\Gamma_{cyl}$ and the inflow and outflow velocity is given by $\bu(0,y) = \bu(2.2,y)= (\frac{6}{0.14^2} \sin(\frac{\pi t}{8})(y(0.41-y)),0)^{'},~~0\le y\le 0.41$.
\end{example}
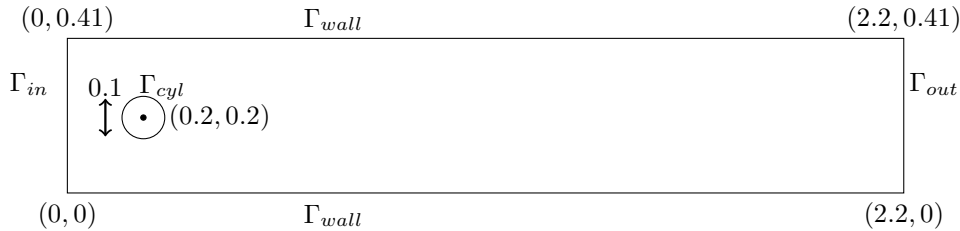
\begin{figure}[!h]
\centering
\begin{tikzpicture}

\draw (0,0) rectangle (5*2.2,5*0.41);
\draw (5*0.2,5*0.2) circle (8pt);
\filldraw (5*0.2,5*0.2) circle (1pt);
 
\node at (0, -0.3) {$(0,0)$};
\node at (5*2.2, -0.3) {$(2.2,0)$};
\node at (0, 5*0.46) {$(0,0.41)$};
\node at (5*2.2, 5*0.46) {$(2.2,0.41)$};
%\node at (5*1.1, -0.3) {$u_1=0,u_2=0$};
%\node at (5*1.1, 5*0.46) {$u_1=0,u_2=0$};
%\node at (-0.7, 5*.20) {$u_1=u_d$};
%\node at (-0.7, 5*.13) {$u_2=0$};
%\node at (5*2.32, 5*.20) {$u_1=u_d$};
%\node at (5*2.32, 5*.13) {$u_2=0$};
%\node at (1, 0.6) {$u_1=0$};
%\node at (1, 0.3) {$u_2=0$};
\node at (5*.4, 1) {$(0.2,0.2)$};
\draw [<->,line width=0.3mm] (5*0.1,5*0.15) -- (5*0.1,5*0.25);
\node at (5*0.1,5*0.28) {$0.1$};

\node at (5*0.7, -0.3) {$\Gamma_{wall}$};
\node at (5*0.7, 5*0.46) {$\Gamma_{wall}$};
\node at (-0.5, 5*.29) {$\Gamma_{in}$};
\node at (5*2.28, 5*.29) {$\Gamma_{out}$};
\node at (5*0.25,5*0.28) {$\Gamma_{cyl}$};

\end{tikzpicture}
\caption{Domain $\Omega$ for flow past cylinder.}
\label{fpc}
\end{figure}

\begin{figure}[!h]
     \centering
     \includegraphics[scale=1.2]{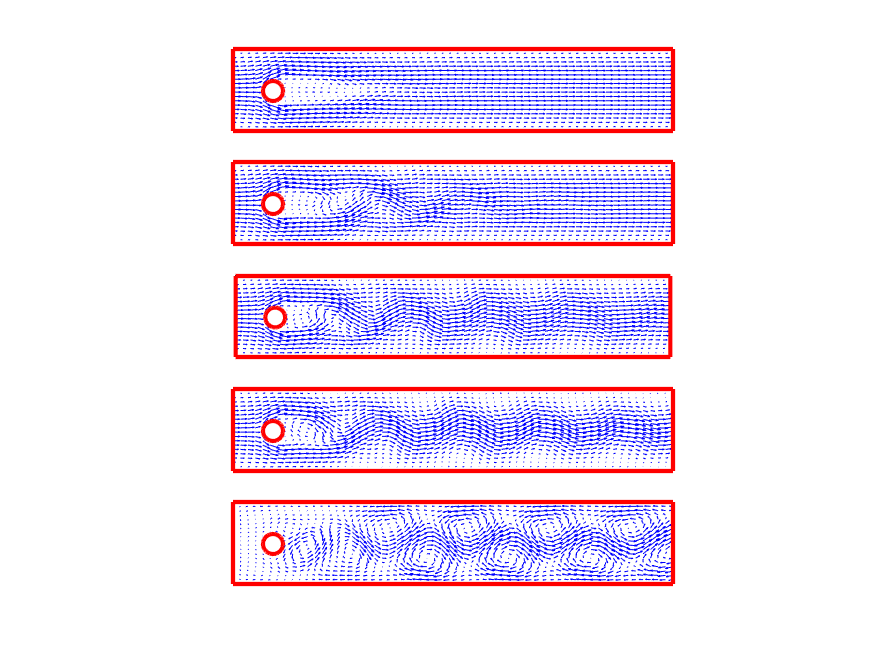} 
     \caption{Velocity field for Example \ref{ex2f} for $T=4,5,6,7,8$.}
     \label{figL3vf}
\end{figure}

\begin{figure}[!h]
     \centering
     \includegraphics[scale=1.2]{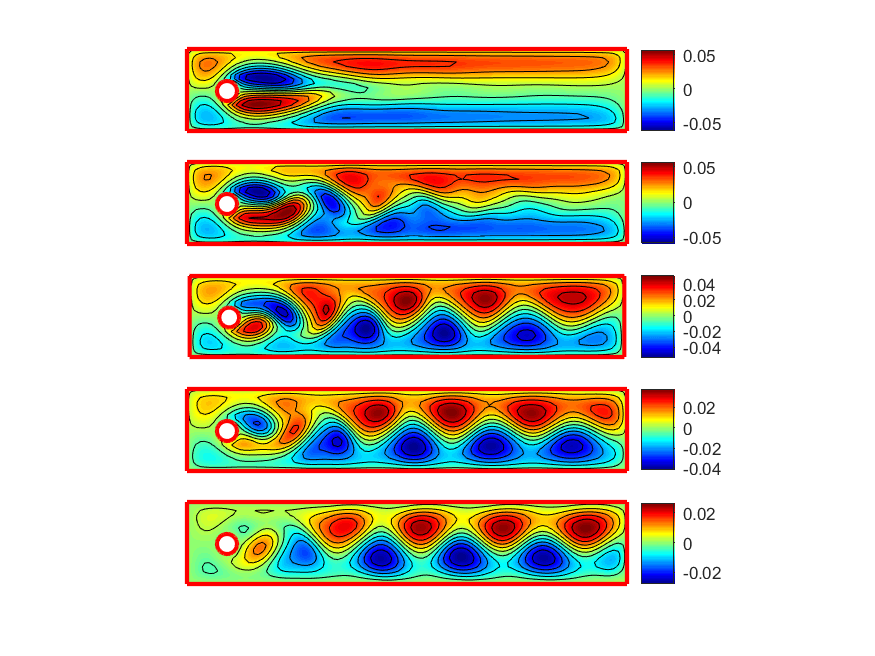} 
     \caption{Stream function for Example \ref{ex2f} for $T=4,5,6,7,8$.}
     \label{figL3vc}
\end{figure}

First, we approximate the velocity and pressure by Taylor-Hood element $P_2-P_1$,  and the domain is discretized with mesh size $h=1/64$. 
For this test, we choose $\nu=10^{-3}$, $k = 10^{-3}$ and the time interval $[0,8]$.
It is known that a vortex sheet develops at the cylinder's bottom around t = 4. In fact, in Figures \ref{figL3vf} and \ref{figL3vc}, we observe this phenomenon, where the velocity field and stream function have been described for different times $T=4,5,6,7$ and $8$. Also, we observe that the vortices separated from the cylinder between the time $T=5$ and $T=6$, and the vortices are still visible at time $T=8$. 

\begin{figure}[!h]
     \centering
     \includegraphics[scale=.35]{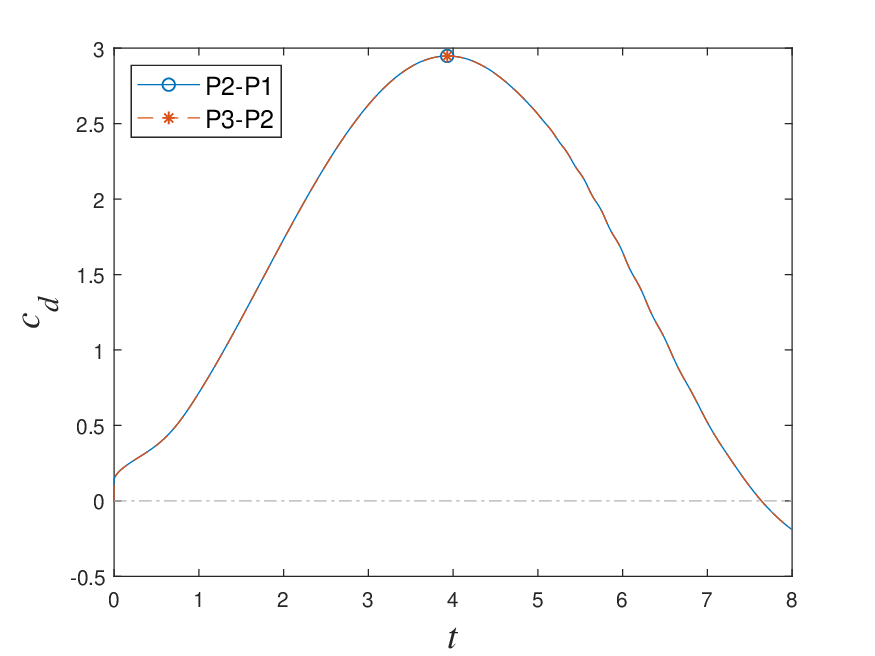} 
     \includegraphics[scale=.35]{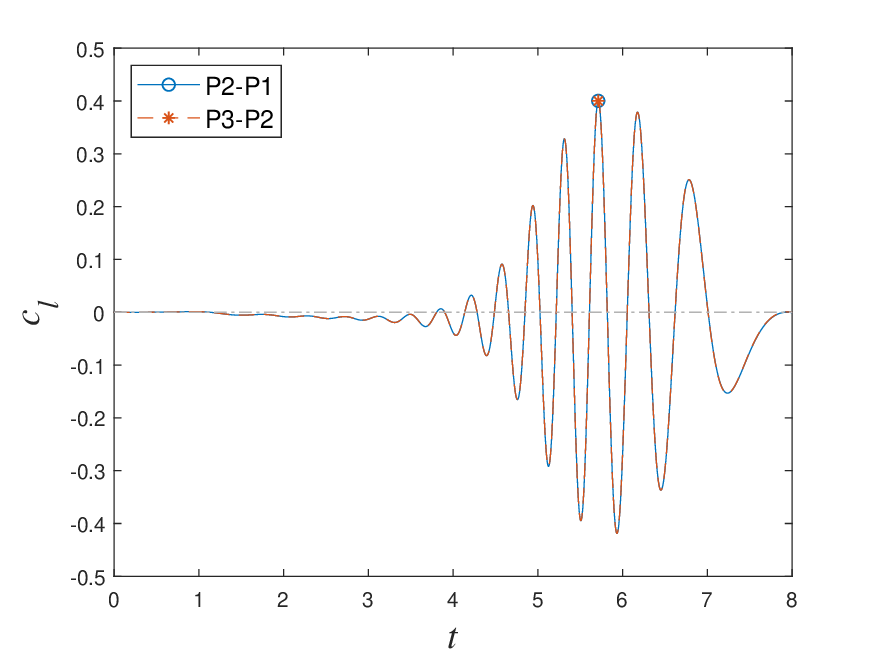} 
     \includegraphics[scale=.35]{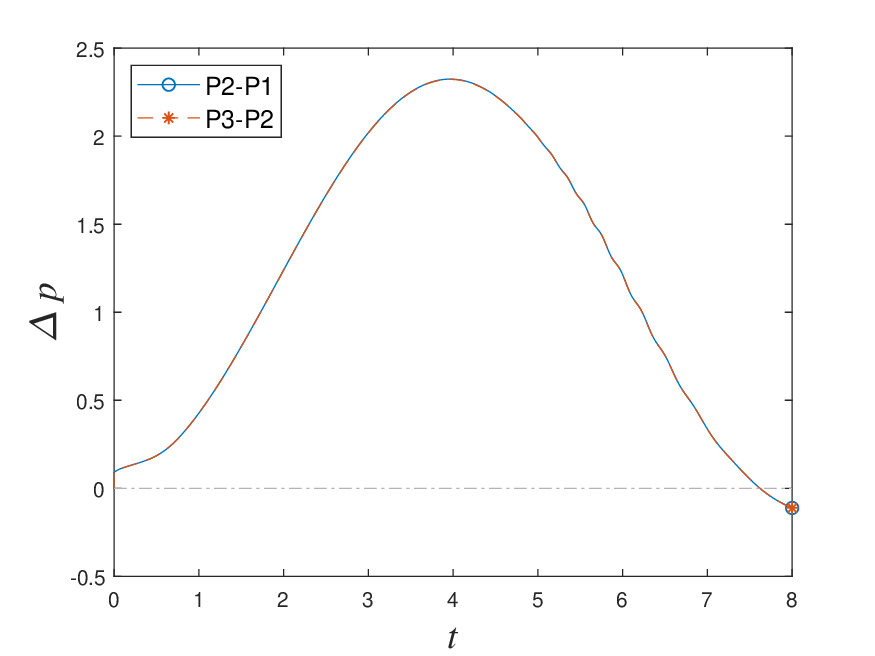} 
     \caption{Drag coefficient, lift coefficient, and pressure difference for Example \ref{ex2f}.}
     \label{figL3v}
\end{figure}

We also plot the evolution of the drag coefficient ($c_d(t)$) at the cylinder, lift coefficient ($c_l(t)$) at the cylinder, and differences in the pressure ($\Delta p(t)$) between the front and the back of the cylinder in Figure \ref{figL3v} for $P_2-P_1$ and $P_3-P_2$ elements. In addition, we mark the maximum value of the drag coefficient, the lift coefficient, and the final value of the pressure difference.  We calculate all these parameters using the formula given in \cite{J04}.


\noindent
In all cases, computation were done in FreeFem++ \cite{H12}.

\section{Conclusion}
This paper deals with  a penalty finite element method along with the backward Euler method  for the incompressible NSEs. With the appropriate use of the inverse of the penalized Stokes operator and the negative norm estimates,  optimal error estimates for the velocity and the pressure terms in both the semi-discrete and fully-discrete schemes are derived. Analysis has been carried out for non-smooth initial data, that is, the initial velocity $\bu_0\in \bH_0^1$. This demands time weighted estimates; the proofs are now more technical and more involved than those  for the smooth case. These optimal estimates are derived under the assumption that the penalized parameter $\ve$ is small. In the numerical part also, optimal convergence rates have been shown for small $\ve.$ Moreover,  results of computational experiments on two benchmark problems show that the proposed method works well for low viscosity and their results compare  well with exiting results from the literature. Although, computational results on one 3D example are encouraging, but this is  with out theoretical justification and this  will form a part of our future endeavour.   

\section*{Acknowledgments}
The first author would like to express his gratitude to the Department of Science and Technology (DST), Government of
India, for the financial support (DST/INSPIRE Fellowship/IF170401).
%
%
%\section*{Declarations}
%
%\textbf{Data Availability Statement}\\
%The datasets generated during and analysed during the current study are available from the corresponding author on reasonable request.
%
%\section*{Conflict of Interest}
%The authors declare that they have no conflict of interest.

%\pagebreak
\section*{Appendix}
\se

\noindent 
\textbf{Proof of the Lemma \ref{pap}:}
\begin{proof} 
Choose $\bphi=\bue$ in (\ref{wfpp}) and use the Cauchy-Schwarz inequality and the Poincar\'e inequality with Lemma \ref{Aep} ($\|\bue\|^2 \le \frac{1}{\lambda_1}\|\nabla\bue\|^2\le \frac{c_0^2}{\lambda_1}\|A_{\ve}^{\frac{1}{2}}\bue\|^2$) to find that
\begin{align}\label{pap000p}
 \frac{d}{dt}\|\bue\|^2  + \nu \|A_{\ve}^{\frac{1}{2}}\bue\|^2  \le  \frac{c_0^2}{\nu\lambda_1}\|\f\|^2.
\end{align}
Note that the non-linear term vanishes due to (\ref{tbp1}). Now multiply by $e^{2\alpha t}$ and integrate from $0$ to $t$ to obtain
\begin{align}\label{pap001p}
 e^{2\alpha t}\|\bue(t)\|^2+(\nu-\frac{2c_0^2\alpha}{\lambda_1})\int_0^t e^{2\alpha s}\|A_{\ve}^{\frac{1}{2}}\bue(s)\|^2ds \le 
\|\bu_{\ve 0}\|^2+\frac{c_0^2(e^{2\alpha t}-1)}{2\alpha\nu\lambda_1}\|\f\|^2_{\infty}.
\end{align}
With $0< \alpha < \frac{\nu\lambda_1}{2c_0^2}$, we have $(\nu-\frac{2c_0^2\alpha}{\lambda_1})>0$. Multiply through out by $e^{-2\alpha t}$ to conclude the first proof. 
Now, we integrate (\ref{pap000p}) with respect to time from $t$ to $t+T$ for any $T>0$, we have
\begin{align}\label{pap0015p}
\|\bue(t+T)\|^2+\nu\int_t^{t+T} \|A_{\ve}^{\frac{1}{2}}\bue(s)\|^2ds \le 
\|\bue(t)\|^2 + \frac{c_0^2T}{\nu\lambda_1}\|\f\|^2_{\infty}.
\end{align}
For the second estimate, choose $\bphi=A_{\ve}^{m+1}\bue$ in (\ref{wfpp}). When $m=0$, we find that
\begin{align}\label{pap01p}
 \frac{1}{2}\frac{d}{dt} \|A_{\ve}^{\frac{1}{2}}\bue\|^2 +\nu\|A_{\ve}\bue\|^2 = (\f,A_{\ve}\bue)  -\tb(\bue,\bue,A_{\ve}\bue).
\end{align}
A use of Ladyzhenskaya's inequality \cite{T84}
($\|\bphi\|_{L^4} \le C \|\bphi\|^{\frac{1}{2}}\|\nabla\bphi\|^{\frac{1}{2}}$, and $\|\nabla\bphi\|_{L^4} \le C \|\nabla\bphi\|^{\frac{1}{2}}\|\Delta\bphi\|^{\frac{1}{2}}$)
with Lemma \ref{Aep}, the Young's inequalities, we bound the nonlinear term as
\begin{align} \label{pap011p}
 \tb(\bue,\bue,A_{\ve}\bue) 
&\le \|\bue\|_{L^4}\|\nabla\bue\|_{L^4} \|A_{\ve}\bue\|  \nonumber\\
&\le C  \|\bue\|^{\frac{1}{2}}\|A_{\ve}^{\frac{1}{2}}\bue\|\|A_{\ve}\bue\|^{3/2} \nonumber\\
&\le C \|\bue\|^2\|A_{\ve}^{\frac{1}{2}}\bue\|^4 + \frac{\nu}{4}\|A_{\ve}\bue\|^2. 
\end{align}
Substitute the above estimate in (\ref{pap01p}) to find that
\begin{align}\label{pap02p}
 \frac{d}{dt} (\|A_{\ve}^{\frac{1}{2}}\bue\|^2) +\nu \|A_{\ve}\bue\|^2 
 \le  C  \big(\|\bue\|^2\|A_{\ve}^{\frac{1}{2}}\bue\|^2\big) \|A_{\ve}^{\frac{1}{2}}\bue\|^2 + \frac{2}{\nu}\|\f\|^2).
\end{align}
We now apply uniform Gronwall's Lemma (Lemma \ref{ugl}) in (\ref{pap02p})  and use (\ref{pap001p}) and (\ref{pap0015p}) to conclude that $\|A_{\ve}^{\frac{1}{2}}\bue(t)\|^2$ is uniformly bounded with respect to $t$ on $[T,\infty)$. Precisely
\begin{equation}\label{pap4p}
  \|A_{\ve}^{\frac{1}{2}}\bue(t)\|^2  \le C, \quad \forall t\ge T.
\end{equation}
For $0\le t\le T$, we use the classical Gronwall's lemma \cite{HR90,S90} in (\ref{pap02p}) and obtain 
\begin{equation}\label{pap41p}
  \|A_{\ve}^{\frac{1}{2}}\bue(t)\|^2  \le C, \quad \text{for}~~ 0\le t\le T.
\end{equation}
Finally,  multiply (\ref{pap02p}) by $e^{2\alpha t}$ and integrate with respect to time from $0$ to $t$ and use the estimates (\ref{pap001p}), (\ref{pap4p}) and (\ref{pap41p}) to complete the second proof when $r=0$. 
For $r=1$, 
we need some intermediate estimate. First we take $\bphi=e^{2\alpha t}\buet$ with $\hbue=e^{\alpha t}\bue$ in (\ref{wfpp}) to obtain
\begin{align}\label{ppap01p}
 \frac{\nu}{2}\frac{d}{dt} \|A_{\ve}^{\frac{1}{2}}\hbue\|^2+  \|\hbuet\|^2= \alpha \nu\|A_{\ve h}^{\frac{1}{2}}\hbue\|^2 
  + (\hf,\buet) -e^{2\alpha t}\tb(\bue,\bue,\buet).
\end{align}
We can estimate the nonlinear term on the right hand side of (\ref{ppap01p}) similar to (\ref{pap011p}) and 
integrate both sides with respect to time to find that
\begin{align*}%\label{ppap02p}
 \nu \|A_{\ve}^{\frac{1}{2}}\hbue\|^2 + \int_0^t\|\hbuet\|^2 ds
  \le C \bigg[\int_0^t \Big(\|A_{\ve}^{\frac{1}{2}}\hbue\|^2 + \|\hf\|^2+ \|A_{\ve}^{\frac{1}{2}}\bue\|^2\|A_{\ve}\hbue\|^2 \big) ds\bigg].
\end{align*}
Now a use of (\ref{pap001p}) and (\ref{pap4p}) lead us to the intermediate estimate. 
\begin{align} \label{ppap021p}
  \|A_{\ve}^{\frac{1}{2}}\bue(t)\|^2 + e^{-2\alpha t}\int_0^t e^{2\alpha s}\|\bues(s)\|^2 ds \le C.
\end{align}
We now differentiate (\ref{wfpp}) with respect to time and deduce that
\begin{align}\label{dwfppp}
 (\buett,\bphi)+\nu a_{\ve}(\buet,\bphi) = (\f_t,\bphi) -\tb(\buet,\bue,\bphi)-\tb(\bue,\buet,\bphi), \quad \forall \bphi\in \bH_0^1.
\end{align}
Take $\bphi=\sigma(t)\buett$ in (\ref{dwfppp}) and use Lemma \ref{Aeph}, the Cauchy-Schwarz inequality to reach at 
\begin{align*}
\frac{d}{dt}(\sigma(t)\|\buet\|^2) + \nu \sigma(t)\|A_{\ve}^{\frac{1}{2}}\buet\|^2 \le C e^{2\alpha t}\|\buet\|^2 +C\sigma(t)\Big(\|\f_t\|^2+ \|\buet\|^2\|A_{\ve h}^{\frac{1}{2}}\bue\|^2\Big).
\end{align*}
Integrate  with respect to time and use (\ref{ppap021p}), (\ref{pap4p}) and (\ref{pap41p}) to obtain
\begin{align} \label{ppap022p}
  \tau(t)\|\buet(t)\|^2 + \nu  e^{-2\alpha t}\int_0^t \sigma(s)\|A_{\ve}^{\frac{1}{2}}\bues(s)\|^2 ds \le C.
\end{align}
Now we are in position to complete the proof of the second estimate when $m=1$. For this, we set $\bphi=A_{\ve}\bue$ in (\ref{wfpp}) and rewrite it and use (\ref{pap011p}) and the Cauchy-Schwarz inequality to arrive at
\begin{align*}
\nu \|A_{\ve h}\bue\|^2 &= (\f,A_{\ve}\bue)-(\buet,A_{\ve}\bue)-\tb(\bue,\bue,A_{\ve}\bue)\\
 &\le  C  \big( \|\f\|^2 +\|\buet\|^2  + \|\bue\|^2\|A_{\ve}^{\frac{1}{2}}\bue\|^4\big) + \frac{\nu}{2}\|A_{\ve}\bue\|^2.
\end{align*}
Multiply by $\tau(t)$ and use (\ref{pap001p}), (\ref{pap4p}), (\ref{pap41p}) and (\ref{ppap022p}) to complete the second proof.
\end{proof}

\noindent
\textbf{Proof of the well-posedness of the discrete solution of problem (\ref{fdbej}):}

\begin{proof}
We can rewritten (\ref{fdbej}) as
\begin{equation*}%\label{wp1}
 (\bUe^n,\bphi_h)+\nu k a_{\ve}(\bUe^n,\bphi_h)+k\tb(\bUe^n,\bUe^n,\bphi_h)= k(\f^n,\bphi_h)
 +(\bUe^{n-1},\bphi_h),~~~\forall\bphi_h\in\bH_h.
\end{equation*}
Consider a function $F:\bH_h\to\bH_h$ such that 
\begin{align*}
(F(\bv),\bphi_h) = (\bv,\bphi_h)+\nu k a_{\ve}(\bv,\bphi_h)+k\tb(\bv,\bv,\bphi_h)- k(\f^n,\bphi_h)
 -(\bUe^{n-1},\bphi_h),~~~\forall\bphi_h\in\bH_h.
\end{align*}
Clearly, $F$ is continuous. Then, a use of pointcar\'e inequality and inverse hypothesis yields
\begin{align*}
(F(\bUe^n),\bUe^n)
&= (\bUe^n,\bUe^n)+\nu k a_{\ve}(\bUe^n,\bUe^n)+k\tb(\bUe^n,\bUe^n,\bUe^n)- k(\f^n,\bUe^n)
 -(\bUe^{n-1},\bUe^n)\\
 &\ge \|\bUe^n\|^2+\nu k \|A_{\ve h}^{\frac{1}{2}} \bUe^n\|^2 - k\|\f^n\|\|\bUe^n\| -\|\bUe^{n-1}\|\|\bUe^n\|\\
 &\ge \Big(\big(1+\frac{\nu k\lambda_1}{c^2}\big) \|\bUe^n\| - \big(k\|\f^n\|+\|\bUe^{n-1}\|\big)\Big)\|\bUe^n\|.
\end{align*}
Now, choose $\bUe^n\in \bH_h$ such that
 $$\bUe^n = \frac{2\big(k\|\f^n\|+\|\bUe^{n-1}\|\big)}{(\big(1+\frac{\nu k\lambda_1}{c^2}\big)}= \alpha_1.$$
 If either $\|\f^n\|\ne 0$ or $\|\bUe^{n-1}\|\ne 0$, then $\alpha_1>0$, which implies that there exists $\bUe^*\in \bH_h$ such that $\|\bUe^*\|\le \alpha_1$ and $F(\bUe^*)=0$.
\end{proof}

\end{document}